\documentclass[11pt]{article}
\usepackage{amsfonts}
\usepackage{color}
\newtheorem{theo}{Theorem}[section]
\newtheorem{lem}[theo]{Lemma}

\newtheorem{defi}[theo]{Definition}

\newcommand{\mysection}[1]{\section{#1} \setcounter{equation}{0}}
\newcommand{\proof}{{\sc Proof.} \quad}
\newcommand{\proofc}{{\sc Proof} \ }
\newcommand{\be}{\begin{equation} \label}
\newcommand{\ee}{\end{equation}}
\newcommand{\bea}{\begin{eqnarray}\label}
\newcommand{\eea}{\end{eqnarray}}
\newcommand{\bas}{\begin{eqnarray*}}
\newcommand{\eas}{\end{eqnarray*}}
\newcommand{\bit}{\begin{itemize}}
\newcommand{\eit}{\end{itemize}}
\newcommand{\qed}{\hfill$\Box$ \vskip.2cm}
\newcommand{\nn}{\nonumber}
\newcommand{\R}{\mathbb{R}}
\newcommand{\N}{\mathbb{N}}
\newcommand{\pO}{\partial\Omega}

\newcommand{\eps}{\varepsilon}

\newcommand{\wto}{\rightharpoonup}
\newcommand{\wsto}{\stackrel{\star}{\rightharpoonup}}
\newcommand{\hra}{\hookrightarrow}
\newcommand{\io}{\int_\Omega}

\newcommand{\mult}{\otimes}

\newcommand{\abs}{\\[5pt]}
\newcommand{\Abs}{\\[5mm]}

\newcommand{\proj}{{\cal P}}
\newcommand{\neps}{n_\eps}
\newcommand{\ceps}{c_\eps}
\newcommand{\ueps}{u_\eps}
\newcommand{\veps}{v_\eps}
\newcommand{\Peps}{P_\eps}
\newcommand{\yeps}{Y_\eps}
\newcommand{\tme}{T_{max,\eps}}

\newcommand{\one}{{\bf{1}}}
\newcommand{\F}{{\cal F}_\kappa}
\newcommand{\cgp}{C_g^+}
\newcommand{\cgm}{C_g^-}
\newcommand{\kz}{K_0}

\begin{document}
\enlargethispage{10mm}
\title{Global weak solutions in a three-dimensional chemotaxis-Navier-Stokes system}
\author{
Michael Winkler\footnote{michael.winkler@math.uni-paderborn.de}\\
{\small Institut f\"ur Mathematik, Universit\"at Paderborn,}\\
{\small 33098 Paderborn, Germany} }
\date{}
\maketitle
\begin{abstract}
\noindent 
  The chemotaxis-Navier-Stokes system 
    \be{00}
    \left\{ \begin{array}{rcll}
    n_t + u\cdot\nabla n &=& \Delta n - \nabla \cdot (n\chi(c)\nabla c),\\[1mm]
    c_t + u\cdot\nabla c &=& \Delta c-nf(c),    \\[1mm]
    u_t + (u\cdot\nabla)u &=& \Delta u +\nabla P + n \nabla \Phi,  \\[1mm]
    \nabla \cdot u &=& 0,
    \end{array} \right.
	\qquad \qquad (\star)
  \ee
  is considered under homogeneous boundary conditions of Neumann type for $n$ and $c$, and of Dirichlet type for $u$,
  in a bounded convex domain $\Omega\subset \R^3$ with smooth boundary, where $\Phi\in W^{1,\infty}(\Omega)$,
  and where $f\in C^1([0,\infty))$ and $\chi\in C^2([0,\infty))$ are nonnegative with $f(0)=0$.
  Problems of this type have been used to describe the mutual interaction of populations of swimming aerobic bacteria
  with the surrounding fluid. Up to now, however, global existence results seem to be available only
  for certain simplified variants such as e.g.~the two-dimensional analogue of ($\star$), or the associated chemotaxis-Stokes
  system obtained on neglecting the nonlinear convective term in the fluid equation.\abs
  The present work gives an affirmative answer to the question of global solvability for ($\star$) in the following sense:
  Under mild assumptions on the initial data, and under modest
  structural assumptions on $f$ and $\chi$, inter alia allowing for the prototypical case when
  \bas
	f(s)=s
	\quad \mbox{for all } s\ge 0
	\qquad \mbox{and} \qquad
	\chi \equiv const.,
  \eas
  the corresponding initial-boundary value problem is shown to possess a globally defined weak solution.\abs
  This solution is obtained as the limit of smooth solutions to suitably regularized problems,
  where appropriate compactness properties are derived
  on the basis of a priori estimates gained from an 		
  energy-type inequality for ($\star$) which in an apparently novel manner combines
  the standard $L^2$ dissipation property of the fluid evolution
  with a quasi-dissipative structure associated with the chemotaxis subsystem in ($\star$).\abs
\noindent {\bf Key words:} chemotaxis, Navier-Stokes, global existence\\
 {\bf AMS Classification:} 35Q92 (primary); 35A01, 35D30,  35Q30 (secondary)
\end{abstract}
\newpage
\section{Introduction}\label{intro}
We consider the chemotaxis-Navier-Stokes system
\be{0}
    \left\{ \begin{array}{rcll}
    n_t + u\cdot\nabla n &=& \Delta n - \nabla \cdot (n\chi(c)\nabla c),
	\qquad & x\in\Omega, \ t>0,\\[1mm]
    c_t + u\cdot\nabla c &=& \Delta c-nf(c), \qquad & x\in\Omega, \ t>0,   \\[1mm]
     u_t + (u\cdot\nabla )u  &=& \Delta u + \nabla P + n \nabla \Phi, \qquad & x\in\Omega, \ t>0, \\[1mm]
    \nabla \cdot u &=& 0, \qquad & x\in\Omega, \ t>0,
    \end{array} \right.
\ee
in a domain $\Omega \subset \R^N$, where the main focus of this work will be on the case when $N=3$ and $\Omega$ is bounded
and convex with smooth boundary.\abs
As described in \cite{goldstein2004}, problems of this type arise in the modeling of populations of swimming
aerobic bacteria in situations
when besides their chemotactically biased movement toward oxygen as their nutrient,
a buoyancy-driven effect of bacterial mass on the fluid motion is not negligible.
Indeed, striking experimental findings indicate that such a mutual chemotaxis-fluid interaction may lead to 
quite complex types of collective behavior, even in markedly simple settings such as present when
populations of {\em Bacillus subtilis} are suspended in 
sessile drops of water (\cite{goldstein2004}, \cite{goldstein2005}, \cite{lorz_M3AS}).\abs
In particular, in (\ref{0}) it is assumed that the presence of bacteria, with density denoted by $n=n(x,t)$,
affects the fluid motion, as represented by its velocity field
$u=u(x,t)$ and the associated pressure $P=P(x,t)$, through buoyant forces.
Moreover, it is assumed that both cells and oxygen, the latter with concentration $c=c(x,t)$, 
are transported by the fluid and diffuse randomly, that the cells partially direct their movement
toward increasing concentrations of oxygen, and that the latter is consumed by the cells.\abs
{\bf The regularity problem in the Navier-Stokes and chemotaxis subsystems.} \quad
The mathematical understanding of such types of interplay is yet quite rudimentary only,
which may be viewed as reflecting the circumstance that (\ref{0}) joins two delicate subsystems which 
themselves 
are far from understood even when decoupled from each other: 
Indeed, as is well-known, 	
the three-dimensional Navier-Stokes system is still lacking a satisfactory existence theory even in absence of external 
forcing terms (\cite{wiegner}):
Global weak solutions for initial data in $L^2(\Omega)$ have been known to exist since Leray's celebrated pioneering 
work (\cite{leray}, \cite{sohr}), 
but despite intense research over the past decades it cannot be decided up to now	
whether the nonlinear convective term may enforce the spontaneous emergence of singularities in the sense
of blow-up with respect to e.g.~the norm in $L^\infty(\Omega)$,
or whether such phenomena are entirely ruled out by diffusion; 
in contrast to this, the latter is known to be the case in the two-dimensional analogue in which unique global 
smooth solutions exist for all reasonably regular initial data to the corresponding Dirichlet problem in bounded domains,
for instance (\cite{sohr}).\abs
A similar criticality of the spatially three-dimensional setting
with respect to rigorous analytical evi\-dence
can be observed for the chemotaxis subsystem obtained upon neglecting the fluid interaction in (\ref{0}).
In fact, e.g.~for the prototypical system
\be{ct}
    \left\{ \begin{array}{rcll}
    n_t &=& \Delta n - \nabla \cdot (n\nabla c),
	\qquad & x\in\Omega, \ t>0,\\[1mm]
    c_t  &=& \Delta c-nc, \qquad & x\in\Omega, \ t>0, 
    \end{array} \right.
\ee
it is known that the Neumann initial-boundary value problem in planar bounded convex domains is uniquely globally solvable
for all suitably smooth nonnegative initial data, whereas in the three-dimensional counterpart only 
certain weak solutions are known to exist globally, with the question whether or not blow-up may occur being
undecided yet (\cite{taowin_JDE2012}).
Anyhow, a highly destabilizing potential of cross-diffusive terms of the type in (\ref{ct}), at relative strength
increasing with the spatial dimension, is indicated by known results on the related classical Keller-Segel system
of chemotaxis, as obtained by replacing the second equation in (\ref{ct}) with $c_t=\Delta c - c + n$:
While all classical solutions to the corresponding initial-boundary value problem remain bounded when either $N=1$, or $N=2$
and the total mass $\io n_0$ of cells is small (\cite{osaki_yagi}, \cite{nagai_senba_yoshida}), 
it is known that finite-time blow-up does occur for large classes of radially symmetric initial data when either
$N=2$ and $\io n_0$ is large, or $N\ge 3$ and $\io n_0$ is an arbitrarily small prescribed number (\cite{mizoguchi_win},
\cite{win_JMPA}).\abs
{\bf Existence results for chemotaxis-fluid systems.} \quad
Accordingly, the literature on coupled chemo\-taxis-fluid systems is yet quite fragmentary,
and most results available so far concentrate either on special cases involving somewhat restrictive assumptions,
or on variants of (\ref{0}) which contain additional regularizing effects.
For instance, a considerable simplification consists in		
removing the convective term $(u\cdot\nabla)u$ from the third equation in (\ref{0}),
thus assuming the fluid motion to be governed by the linear Stokes equations.
The correspondingly modified system is indeed known to possess global solutions at least in a certain weak
sense under suitable initial and boundary conditions in smoothly bounded three-dimensional convex domains,
provided that the coefficient functions in (\ref{0}) are adequately smooth, and that
$\chi$ and $f$ satisfy some mild structural conditions (cf.~(\ref{struct}) below) generalizing the prototypical choices
\be{proto}
	f(s)=s \quad \mbox{for all } s\ge 0
	\qquad \mbox{and} \qquad
	\chi\equiv const.
\ee
(see (\cite{win_CPDE});
it is not known, however, whether these solutions are sufficiently regular so as to avoid phenomena of unboundedness,
e.g.~with respect to the norm of $n$ in $L^\infty(\Omega)$, in either finite or infinite time
(cf.~also \cite{chae_kang_lee_cpde} for some refined extensibility criteria for local-in-time smooth solutions).\abs
A further regularization can be achieved by assuming the diffusion of cells to be nonlinearly enhanced at large densities.
Indeed, if in the first equation the term $\Delta n$ is replaced by $\Delta n^m$ for $m>1$, then, firstly, 
for any such choice of $m$ it is again possible to construct global weak solutions under appropriate assumptions on
$\chi$ and $f$ (\cite{duan_xiang_IMRN2012}), 
but beyond this one can secondly prove local-in-time and even global boundedness of these
solutions in the cases $m>\frac{8}{7}$ and $m>\frac{7}{6}$, respectively, and thereby rule out the occurrence of 
blow-up in finite time, and also in infinite time (cf.~\cite{taowin_ANNIHP} and \cite{win_ctf_3d_nonlinear_general}).
If the range of $m$ is further restricited by assuming $m>\frac{4}{3}$, then global existence of, possibly unbounded,
solutions can be derived even in presence of the full nonlinear Navier-Stokes equations (\cite{vorotnikov}).\\
As alternative blow-up preventing mechanisms, the authors in \cite{cao_wang} and \cite{vorotnikov} identify 
certain saturation effects at large cell densities in the cross-diffusive term, as well as the inclusion
of logistic-type cell kinetics with quadratic death terms in (\ref{0}), in both cases leading to corresponding
results on global existence of weak solutions.\abs
In the spatially two-dimensional case, the knowledge on systems of type (\ref{0}) is expectedly much further developed.
Even in the original chemotaxis-Navier Stokes system (\ref{0}) containing nonlinear convection in the fluid evolution, 
the regularizing effect of the diffusive mechanisms turns out to be strong enough so as to allow
for the construction of unique global bounded classical solutions under the mild assumptions (\ref{reg_coeff})
and (\ref{struct}) on $\chi, f$ and $\Phi$ (\cite{win_CPDE}; see also \cite{liu_lorz}),
and to furthermore enforce stabilization of these solutions toward spatially homogeneous
equilibria in the large time limit (\cite{win_ARMA}).
Corresponding results on global existence in presence of porous medium type cell diffusion, or of additional logistic terms,
can be found in \cite{DiFLM}, \cite{taowin_DCDSA}, \cite{vorotnikov}, \cite{cao_ishida_NON} and \cite{ishida}, for instance,
and recently statements on global existence and boundedness 
have been derived in \cite{espejo_suzuki} and \cite{tao_ctf} for a two-dimensional chemotaxis-Stokes
variant of (\ref{0}) involving signal production by cells and a quadratic death term in the cell evolution,
as proposed in a different modeling context in \cite{kiselev_ryzhik_CPDE}.\abs
{\bf Main results.} \quad
For the full three-dimensional chemotaxis-Navier-Stokes system (\ref{0}), 
even at the very basic level of global existence in generalized solution frameworks, a satisfactory solution theory
is entirely lacking.
The only global existence results we are aware of
concentrate on the construction of solutions near constant steady states (\cite{DLM}), or on the 
particular case when $\chi$ precisely coincides with a multiple of $f$ (\cite{chae_kang_lee}), 
where the latter not only excludes the situation determined by (\ref{proto}), but
under the natural assumption that $f(0)=0$ apparently also rules out any choice of $\chi$ which is consistent with
standard approaches in the modeling of chemotaxis phenomena (\cite{hillen_painter2009}).\abs
It is the purpose of the present work to undertake a first step toward a comprehensive existence theory for (\ref{0})
under mild assumptions of the coefficient functions therein, and for widely general initial data.
In order to formulate our main results in this direction, let us specify the precise evolution problem addressed in the sequel
by considering (\ref{0}) along with the initial conditions
\be{0i}
    n(x,0)=n_0(x), \quad c(x,0)=c_0(x) \quad \mbox{and} \quad u(x,0)=u_0(x), \qquad x\in\Omega,
\ee
and under the boundary conditions
\be{0b}
	\frac{\partial n}{\partial\nu}=\frac{\partial c}{\partial\nu}=0 
	\quad \mbox{and} \quad u=0
	\qquad \mbox{on } \pO,
\ee
in a bounded convex domain $\Omega\subset \R^3$ with smooth boundary,
where we assume that 
\be{init}
    \left\{
    \begin{array}{l}
    n_0 \in L \log L(\Omega) \quad \mbox{is nonnegative with $n_0\not\equiv 0$,} \quad
	\mbox{that} \\
    c_0 \in L^\infty(\Omega) \quad \mbox{ is nonnegative and such that $\sqrt{c_0} \in W^{1,2}(\Omega)$, \quad and that}\\
    u_0 \in L^2_\sigma(\Omega),
    \end{array}
    \right.
\ee
with $L^2_\sigma(\Omega):=\{ \varphi\in L^2(\Omega) \ | \ \nabla \cdot \varphi=0 \}$ denoting the Hilbert space
of all solenoidal vector fields in $L^2(\Omega)$.\\
With regard to the chemotacitic sensitivity $\chi$, the signal consumption rate $f$ and the 	
potential $\Phi$ in (\ref{0}), throughout this paper we shall require that
\be{reg_coeff}
	\left\{ \begin{array}{l}
	\chi \in C^2([0,\infty)) \quad \mbox{is positive on $[0,\infty)$}, \\
	f\in C^1([0,\infty)) \quad \mbox{is nonnegative on $[0,\infty)$ with $f(0)=0$, and that} \\
	\Phi\in W^{1,\infty}(\Omega),
	\end{array} \right.
\ee
and moreover we will need the structural hypotheses
\be{struct}
	\Big(\frac{f}{\chi}\Big)'>0,
	\quad 
	\Big(\frac{f}{\chi}\Big)'' \le 0
	\quad \mbox{and} \quad
	(\chi\cdot f)'' \ge 0
	\qquad \mbox{on } [0,\infty).
\ee
We shall see that within this framework, there exists at least one globally defined triple $(n,c,u)$
of functions solving (\ref{0}) in a natural generalized sense specified in Definition \ref{defi_weak} below.
Apart from satisfying the respective weak formulations associated with the PDEs in (\ref{0}), this
solution will enjoy further properties in that it fulfils two energy-type inequalities.
The first of these will be the standard estimate (\ref{energy1}) reflecting
energy dissipation in the Navier-Stokes system, as satisfied by any so-called
turbulent solution thereof (\cite{wiegner}, \cite{sohr}),
while the second will refer to the functional $\F$ with appropriate $\kappa>0$, where we have set
\be{def_energy}
	\F[n,c,u] := \io n\ln n 
	+ \frac{1}{2} \io \frac{\chi(c)}{f(c)} |\nabla c|^2 
	+ \kappa \io |u|^2
\ee
for $\kappa>0$ whenever 
$n\in L\log L(\Omega)$ and $c\in W^{1,2}(\Omega)$ are nonnegative and such that 
$\frac{\chi(c)}{f(c)}|\nabla c|^2 \in L^1(\Omega)$, and $u\in L^2(\Omega;\R^3)$.\abs
Now our main result reads as follows.
\begin{theo}\label{theo_global}
  Let (\ref{reg_coeff}) and (\ref{struct}) hold.
  Then for all $n_0, c_0$ and $u_0$ fulfilling (\ref{init}), there exist 	
  \bea{reg_w}
	& & n\in L^\infty((0,\infty);L^1(\Omega))
	\cap L^\frac{5}{4}_{loc}([0,\infty);W^{1,\frac{5}{4}}(\Omega)), \nn\\
	& & c\in L^\infty(\Omega\times (0,\infty))
	\cap L^4_{loc}([0,\infty);W^{1,4}(\Omega)), \nn\\
	& & u\in L^\infty_{loc}([T,\infty);L^2_\sigma(\Omega))
	\cap L^2_{loc}([T,\infty);W_0^{1,2}(\Omega)),
  \eea
  such that $(n,c,u)$ is a global weak solution of the problem (\ref{0}), (\ref{0i}), (\ref{0b}) 
  in the sense of Definition \ref{defi_weak}.
  This solution can be obtained as the pointwise limit a.e.~in $\Omega\times (0,\infty)$ of a suitable sequence 
  of classical solutions to the regularized problems (\ref{0eps}) below.
  Moreover, $(n,c,u)$ has the additional properties that
  \be{reg_ees}
	\begin{array}{l}
	n^\frac{1}{2} \in L^2_{loc}([0,\infty);W^{1,2}(\Omega)) 
	\quad \mbox{and} \\
	c^\frac{1}{4} \in L^4_{loc}([T,\infty);W^{1,4}(\Omega)),
	\end{array}
  \ee
  and there exist $\kappa>0, K>0$ and a null set $N\subset (0,\infty)$ such that
  \be{energy1}
	\frac{1}{2} \io |u(\cdot,t)|^2 + \int_{t_0}^t \io |\nabla u|^2
	\le \frac{1}{2} \io |u(\cdot,t_0)|^2 
	+\int_{t_0}^t \io nu\cdot\nabla \Phi
	\qquad \mbox{for all $t_0\in [0,\infty)\setminus N$ and all } t>t_0,
  \ee
  as well as
  \bea{energy}
	\frac{d}{dt} \F[n,c,u]
	+ \frac{1}{K} \io \bigg\{ \frac{|\nabla n|^2}{n} + \frac{|\nabla c|^4}{c^3} + |\nabla u|^2 \bigg\}
	\le K
	\qquad \mbox{in ${\mathcal{D}}'([0,\infty))$.}
  \eea
\end{theo}
{\bf Remark.} \quad
  The property $c^\frac{1}{4} \in L^4_{loc}([0,\infty);W^{1,4}(\Omega))$ along with the boundedness
  of $c$ and the fact that $\frac{f}{\chi}$ is nonnegative and belongs to $C^1([0,\infty))$, 
  with nonvanishing derivative at zero,		
  ensures that $\frac{\chi(c)}{f(c)} |\nabla c|^2 \in L^1(\Omega)$ for a.e.~$t>0$ by the Cauchy-Schwarz inequality.
  Along with the regularity features of $n$ and $u$ in (\ref{reg_w}) this implies that
  $\F[n,c,u](t)$ and $\io |u(\cdot,t)|^2$ are is well-defined for a.e.~$t>0$, 
  whence the statements (\ref{energy}) and (\ref{energy1}) are indeed meaningful.\Abs
The plan of this paper is as follows. In Section \ref{sect2} we shall specify the generalized solution concept
considered thereafter, and introduce a family of regularized problems each of which allows for smooth solutions 
at least locally in time.
Section \ref{sect3.1} will be devoted to an analysis of the functional
obtained on letting $\kappa=0$ in (\ref{def_energy}), evaluated at these approximate solutions. 
As known from previous studies, the assumptions in (\ref{struct}) ensure that the time evolution of this 
two-component functional
involves, besides certain dissipated quantities, expressions containing the fluid velocity.
An apparently novel way to treat the latter by making appropriate use of the standard energy dissipation in the Navier-Stokes
equations will allow for absorbing these suitably in Section \ref{sect3.2}.
This will entail a series of a priori estimates which will firstly be used in Section \ref{sect3.3} to make sure
that all the approximate solutions are actually global in time, and which secondly enable us to derive 
further $\eps$-independent bounds in Section \ref{sect3.4}.
On the basis of the compactness properties thereby implied, in Section \ref{sect4} we shall finally pass
to the limit along an adequate sequence of numbers $\eps=\eps_j\searrow 0$ and thereby verify Theorem \ref{theo_global}.
\mysection{Preliminaries}		
\subsection{A weak solution concept}
We first specify the notion of weak solution to which we will refer in the sequel.
Here for candidates of solutions we require the apparently weakest possible regularity properties
which ensure that all expressions in the weak identities (\ref{w1}), (\ref{w2}) and (\ref{w3}) are meaningful.
As already announced in the formulation of Theorem \ref{theo_global}, the solution we shall construct below
will actually be significantly more regular.\\
Throughout the sequel, for vectors $v\in \R^3$ and $w\in\R^3$ we let $v\mult w$ denote the matrix 
$(a_{ij})_{i,j\in \{1,2,3\}}\in\R^{3\times 3}$ defined on setting $a_{ij}:=v_i w_j$ for $i,j\in \{1,2,3\}$.
\begin{defi}\label{defi_weak}
  By a {\em global weak solution} of (\ref{0}), (\ref{0i}), (\ref{0b}) we mean a triple $(n,c,u)$ of functions 
  \bea{reg_weak}
	n\in L^1_{loc}([0,\infty);W^{1,1}(\Omega)), 
	\quad
	c\in L^1_{loc}([0,\infty);W^{1,1}(\Omega)), 
	\quad
	u\in L^1_{loc}([0,\infty);W_0^{1,1}(\Omega;\R^3)), 
  \eea
  such that $n\ge 0$ and $c\ge 0$ a.e.~in $\Omega\times (0,\infty)$,
  \bea{reg_weak2}
	& & nf(c)\in L^1_{loc}(\bar\Omega\times [0,\infty)),
	\qquad 
	u\mult u \in L^1_{loc}(\bar\Omega\times [0,\infty);\R^{3\times 3}),
	\qquad \mbox{and} \nn\\
	& & n\chi(c)\nabla c, nu \mbox{ and } cu
	\ \mbox{belong to } L^1_{loc}(\bar\Omega\times [0,\infty);\R^3),
  \eea
  that $\nabla \cdot u=0$ a.e.~in $\Omega\times (0,\infty)$, and that
  \bea{w1}
	-\int_0^\infty \io n\phi_t - \io n_0\phi(\cdot,0)
	= - \int_0^\infty \io \nabla n\cdot \nabla \phi
	+ \int_0^\infty \io n\chi(c) \nabla c \cdot\nabla\phi
	+ \int_0^\infty \io n u \cdot \nabla\phi
  \eea
  for all $\phi\in C_0^\infty(\bar\Omega\times [0,\infty))$,
  \bea{w2}
	-\int_0^\infty \io c\phi_t - \io c_0\phi(\cdot,0)
	= - \int_0^\infty \io \nabla c\cdot \nabla \phi
	- \int_0^\infty \io nf(c) \phi
	+ \int_0^\infty \io c u\cdot\nabla \phi
  \eea
  for all $\phi\in C_0^\infty(\bar\Omega\times [0,\infty))$ as well as
  \bea{w3}
	-\int_0^\infty \io u\cdot\phi_t - \io u_0\cdot \phi(\cdot,0)
	= - \int_0^\infty \io \nabla u \cdot\nabla \phi
	+ \int_0^\infty u\mult u \cdot \nabla \phi
	+ \int_0^\infty \io n\nabla \Phi\cdot \phi
  \eea
  for all $\phi\in C_0^\infty(\Omega\times [0,\infty);\R^3)$ satisfying $\nabla\cdot \phi\equiv 0$.
\end{defi}
\subsection{A family of regularized problems}\label{sect2}
In order to suitably regularize the original problem (\ref{0}), (\ref{0i}), (\ref{0b}), let us consider families 
of approximate initial data $n_{0\eps}, c_{0\eps}$ and $u_{0\eps}$, $\eps\in (0,1)$, with the properties that
\be{I1}
	\left\{ \begin{array}{l}
	n_{0\eps} \in C_0^\infty(\Omega),
	\quad n_{0\eps}\ge 0 \mbox{ in $\Omega$,}
	\quad
	\io n_{0\eps}=\io n_0
	\quad \mbox{for all } \eps\in (0,1) 
	\qquad \mbox{and} \\[1mm]
	n_{0\eps} \to n_0
	\quad \mbox{in } L\log L(\Omega)
	\qquad \mbox{as } \eps\searrow 0,
	\end{array} \right.
\ee
that
\be{I2}
	\left\{ \begin{array}{l}
	c_{0\eps} \ge 0 \mbox{ in $\Omega$ is such that }
	\sqrt{c_{0\eps}} \in C_0^\infty(\Omega)
	\quad \mbox{and} \quad
	\|c_{0\eps}\|_{L^\infty(\Omega)} \le \|c_0\|_{L^\infty(\Omega)}
	\quad \mbox{for all } \eps\in (0,1) 
	\qquad \mbox{and} \\[1mm]
	\sqrt{c_{0\eps}} \to \sqrt{c_0}
	\quad \mbox{a.e.~in $\Omega$ and in } W^{1,2}(\Omega)
	\qquad \mbox{as } \eps\searrow 0,
	\end{array} \right.
\ee
and that
\be{I3}
	\left\{ \begin{array}{l}
	u_{0\eps} \in C_{0,\sigma}^\infty(\Omega)
	\quad \mbox{with} \quad
	\|u_{0 \eps}\|_{L^2(\Omega)}=\|u_0\|_{L^2(\Omega)}
	\quad \mbox{for all } \eps\in (0,1) 
	\qquad \mbox{and} \\[1mm]
	u_{0\eps} \to u_0
	\quad \mbox{in } L^2(\Omega)
	\qquad \mbox{as } \eps\searrow 0,
	\end{array} \right.
\ee
where as usual $Llog L(\Omega)$ denotes the standard Orlicz space associated with the Young function 
$(0,\infty) \ni z \mapsto z\ln(1+z)$.\\
For $\eps\in (0,1)$, we thereupon consider
\be{0eps}
    \left\{ \begin{array}{rcll}
    n_{\eps t} + \ueps\cdot\nabla \neps &=& \Delta \neps - \nabla \cdot (\neps F'_\eps(\neps)\chi(\ceps)\nabla \ceps),
	\qquad & x\in\Omega, \ t>0,\\[1mm]
    c_{\eps t} + \ueps\cdot\nabla \ceps &=& \Delta \ceps-F_\eps(\neps)f(\ceps), \qquad & x\in\Omega, \ t>0,   \\[1mm]
    u_{\eps t} + (\yeps \ueps\cdot\nabla) \ueps   &=& \Delta \ueps + \nabla \Peps + \neps \nabla \phi, 
	\qquad & x\in\Omega, \ t>0, \\[1mm]
    \nabla \cdot \ueps &=& 0, \qquad & x\in\Omega, \ t>0, \\[1mm]
	& & \hspace*{-32mm}
	\frac{\partial\neps}{\partial\nu}=\frac{\partial\ceps}{\partial\nu}=0, \quad \ueps=0,
	\qquad & x\in\pO, \ t>0, \\[1mm]
	& & \hspace*{-32mm}
	\neps(x,0)=n_{0\eps}(x), \quad \ceps(x,0)=c_{0 \eps}(x), \quad \ueps(x,0)=u_{0\eps}(x),
	\qquad & x\in\Omega,
    \end{array} \right.
\ee
where we adopt from \cite{win_CPDE} the weakly increasing approximation $F_\eps$ 
of $[0,\infty) \ni s \mapsto s$ determined by
\be{Feps}
	F_\eps(s):=\frac{1}{\eps} \ln \Big(s+\frac{1}{\eps}\Big)
	\qquad \mbox{for } s\ge 0,
\ee
and where we utilize the standard Yosida approximation $\yeps$ (\cite{sohr}, \cite{miyakawa_sohr88}) defined by
\be{yeps}
	\yeps v:=(1+\eps A)^{-1} v
	\qquad \mbox{for } v\in L^2_\sigma(\Omega).
\ee
Here and throughout the sequel, by $A$ we mean the realization of the Stokes operator $-\proj \Delta$
in $L^2_\sigma(\Omega)$, with domain $D(A)=W^{2,2}(\Omega) \cap W^{1,2}_{0,\sigma}(\Omega)$, where
$W^{1,2}_{0,\sigma}(\Omega):=W_0^{1,2}(\Omega) \cap L^2_\sigma(\Omega) \equiv 
\overline{C_{0,\sigma}^\infty(\Omega)}^{\, \|\cdot\|_{W^{1,2}(\Omega)}}$ 
with $C_{0,\sigma}^\infty(\Omega):=C_0^\infty(\Omega) \cap L^2_\sigma(\Omega)$, and where
$\proj$ denotes the Helmholtz projection in $L^2(\Omega)$.
It is well-known that $A$ is self-adjoint and positive due to the fact that $\Omega$ is bounded, and hence
in particular possesses fractional powers $A^\alpha$ for arbitrary $\alpha\in\R$ (\cite[Ch.~III.2]{sohr}).\abs
We remark that in contrast to the case of the pure Navier-Stokes equations without chemotactic coupling,
where global existence of weak solutions can be proved employing the less regularizing operators 
$(1+\eps A^\frac{1}{2})^{-1}$ instead of $\yeps$ (\cite[Ch.~V.2]{sohr}), the use of our stronger regularization 
in (\ref{yeps})
will turn out to be more convenient in the present context, because in conjunction with the properties of $F_\eps$
it will allow for a comparatively simple proof
of global solvability in (\ref{0eps}) due to the fact that $\yeps$ acts as a bounded operator from $L^2(\Omega)$
into $L^\infty(\Omega)$ (cf.~Lemma \ref{lem777}).\abs
Let us furthermore note that our choice of $F_\eps$ ensures that
\be{Feps1}
	0 \le F_\eps'(s)=\frac{1}{1+\eps s} \le 1
	\quad \mbox{and} \quad
	0 \le F_\eps(s) \le s
	\qquad \mbox{for all $s\ge 0$ and $\eps\in (0,1)$,}
\ee
and that
\be{Feps2}
	F'_\eps(s) \nearrow 1
	\quad \mbox{and} \quad
	F_\eps(s) \nearrow s
	\quad \mbox{as $\eps\searrow 0$}
	\qquad \mbox{for all $s\ge 0$.}
\ee
All the above approximate problems admit for local-in-time smooth solutions:
\begin{lem}\label{lem90}
  For each $\eps\in (0,1)$, there exist $\tme\in (0,\infty]$ and uniquely determined functions
  \be{90.1}
	\neps \in C^{2,1}(\bar\Omega\times [0,\tme)),
	\quad
	\ceps \in C^{2,1}(\bar\Omega\times [0,\tme))
	\quad \mbox{and} \quad
	\ueps \in C^{2,1}(\bar\Omega\times [0,\tme);\R^3)
  \ee
  which are such that $\neps>0$ and $\ceps>0$ in $\bar\Omega\times (0,\tme)$,
  and such that with some $P_\eps\in C^{1,0}(\Omega\times (0,\tme))$, the quadruple $(\neps,\ceps,\ueps,P_\eps)$
  solves (\ref{0eps}) classically in $\Omega\times (0,\tme)$.
  Moreover, 
  \bea{ext_crit}
	& & \hspace*{-20mm}
	\mbox{if $\tme<\infty$, then }
	\|\neps(\cdot,t)\|_{L^\infty(\Omega)} + \|\ceps(\cdot,t)\|_{W^{1,q}(\Omega)}
	+ \|A^\alpha\ueps(\cdot,t)\|_{L^2(\Omega)}
	\to \infty
	\quad \mbox{as } t\nearrow \tme \nn\\[1mm]
	& & \hspace*{80mm} \mbox{for all $q>3$ and $\alpha>\frac{3}{4}$}.
  \eea
\end{lem}
\proof
  A proof for this, based e.g.~on the contraction mapping principle and standard regularity theories for the
  heat equation and the Stokes system (\cite{LSU}, \cite{quittner_souplet}, \cite{solonnikov2007}, \cite{sohr}),
  can be copied almost word by word from \cite[Lemma 2.1]{win_CPDE}, where minor 
  modifications, mainly due to the presence of the Yosida approximation operator $\yeps$, may be left to the reader.
\qed
We shall later see in Lemma \ref{lem777} that each of these solutions is in fact global in time. 
This will be a particular consequence of a series of a priori estimates for (\ref{0eps}), 
as the first two of which one may view the following two basic properties.
\begin{lem}\label{lem_basic}
  For any $\eps\in (0,1)$, we have
  \be{mass_eps}
	\io \neps(\cdot,t)=\io n_0
	\qquad \mbox{for all } t\in (0,\tme)
  \ee
  and
  \be{cinfty_eps}
	\|\ceps(\cdot,t)\|_{L^\infty(\Omega)} \le s_0:=\|c_0\|_{L^\infty(\Omega)}
	\qquad \mbox{for all } t\in (0,\tme).
  \ee
\end{lem}
\proof
  In (\ref{0eps}), we only need to integrate the first equation over $\Omega$ and apply the maximum principle to 
  the second equation.
\qed
\mysection{A priori estimates}		
\subsection{An energy functional for the chemotaxis subsystem}\label{sect3.1}
A key role in the derivation of further estimates will be played by the following identity which was stated
in \cite{win_CPDE} for the case when $\yeps$ in (\ref{0eps}) is replaced by the identity operator, but which readily
extends to the present situation, because it actually only relies on the first two equations in (\ref{0eps})
and the fact that the fluid component in the transport terms therein is solenoidal.
The novelty of the present reasoning, as compared to previous approaches based on this or similar identities
(cf.~also \cite{DLM}), appears to consist in the particular manner in which (\ref{72.1}) will subsequently 
be related to the natural dissipative properties of the Navier-Stokes subsystem in (\ref{0eps}).
\begin{lem}\label{lem72}
  Given any $\eps\in (0,1)$, the solution of (\ref{0eps}) satisfies
  \bea{72.1}
	& & \hspace*{-20mm}
	\frac{d}{dt} \bigg\{ \io \neps \ln \neps + \frac{1}{2} \io |\nabla \Psi(\ceps)|^2 \bigg\}
	+ \io \frac{|\nabla\neps|^2}{\neps} 
	+ \io g(\ceps)|D^2\rho(\ceps)|^2 \nn\\
	&=& - \frac{1}{2} \io \frac{g'(\ceps)}{g^2(\ceps)} |\nabla \ceps|^2 (\ueps\cdot\nabla \ceps)
	+ \io \frac{1}{g(\ceps)} \Delta\ceps (\ueps\cdot\nabla\ceps) \nn\\
	& & + \io F_\eps(\neps) \Big\{ \frac{f(\ceps) g'(\ceps)}{2g^2(\ceps)} - \frac{f'(\ceps)}{g(\ceps)} \Big\}
	\cdot |\nabla\ceps|^2 \nn\\
	& & + \frac{1}{2} \io \frac{g''(\ceps)}{g^2(\ceps)} |\nabla \ceps|^4
	+ \frac{1}{2} \int_{\pO} \frac{1}{g(\ceps)} \cdot \frac{\partial |\nabla\ceps|^2}{\partial\nu}
	\qquad \mbox{for all } t\in (0,\tme),
  \eea
  where we have set
  \be{g_Psi_rho}
	g(s):=\frac{f(s)}{\chi(s)},
	\quad
	\Psi(s):=\int_1^s \frac{d\sigma}{\sqrt{g(\sigma)}}
	\quad \mbox{and} \quad
	\rho(s):=\int_1^s \frac{d\sigma}{g(\sigma)}
	\qquad \mbox{for } s>0.
  \ee
\end{lem}
\proof
  This can be obtained by straightforward computation on the basis of the first two equations in (\ref{0eps}) and the
  fact that $\nabla\cdot\ueps\equiv 0$ (see \cite[Lemma 3.2]{win_CPDE} for details).
\qed
In order to take full advantage of the dissipated quantities on the left of (\ref{72.1}), we recall the following
functional inequality from \cite[Lemma 3.3]{win_CPDE}.
\begin{lem}\label{lem73}
  Suppose that $h\in C^1((0,\infty))$ is positive, and let $\Theta(s):=\int_1^s \frac{d\sigma}{h(\sigma)}$ for $s>0$.
  Then 
  \be{73.1}
	\io \frac{h'(\varphi)}{h^3(\varphi)} |\nabla \varphi|^4 
	\le (2+\sqrt{3})^2 \io \frac{h(\varphi)}{h'(\varphi)} |D^2 \Theta(\varphi)|^2
  \ee
  holds for all $\varphi\in C^2(\bar\Omega)$ satisfying $\varphi>0$ in $\bar\Omega$ and 
  $\frac{\partial\varphi}{\partial\nu}=0$ on $\pO$.
\end{lem}
For our application of (\ref{73.1}) to (\ref{72.1}), let us state a consequence of our hypotheses (\ref{reg_coeff})
and (\ref{struct}) on the function $g$ from (\ref{g_Psi_rho}) which appears as a weight function in several expressions
in (\ref{72.1}).
\begin{lem}\label{lem71}
  Let $s_0$ be as in (\ref{cinfty_eps}).
  Then there exist $\cgp>0$ and $\cgm>0$ such that the function $g=\frac{f}{\chi}$ in (\ref{g_Psi_rho}) satisfies
  \be{71.1}
	\cgm \cdot s \le g(s) \le \cgp \cdot s
	\qquad \mbox{for all } s\in [0,s_0].
  \ee
\end{lem}
\proof 
  This is an immediate consequence of the first assumption in (\ref{struct}), which together with (\ref{reg_coeff})
  entails that $g$ belongs to $C^1([0,s_0])$ with $g(0)=0, g'(0)>0$ and $g>0$ on $(0,s_0]$.
\qed
We can thereby turn (\ref{72.1}) into an inequality only involving the Dirichlet integral of the fluid velocity
on its right-hand side.
\begin{lem}\label{lem74}
  There exists $\kz\ge 1$ such that for all $\eps\in (0,1)$ we have
  \bea{74.1}
	\frac{d}{dt} \bigg\{ \io \neps \ln \neps + \frac{1}{2} \io |\nabla \Psi(\ceps)|^2 \bigg\}
	&+& \frac{1}{\kz} \cdot \bigg\{ \io \frac{|\nabla\neps|^2}{\neps} 
	+ \io \frac{|D^2 \ceps|^2}{\ceps}
	+ \io \frac{|\nabla \ceps|^4}{\ceps^3} \bigg\} \nn\\[1mm]
	&\le& \kz\io |\nabla\ueps|^2
  \eea
  for all $t\in (0,\tme)$, where $\Psi$ is as in (\ref{g_Psi_rho}).
\end{lem}
\proof 
  We first follow an argument presented in \cite[Lemma 3.4]{win_CPDE} to infer from the third and the second
  inequality in (\ref{struct}) that with $g=\frac{f}{\chi}$ we have
  \be{74.111}
	\frac{fg'}{2g^2} - \frac{f'}{g} = - \frac{(\chi\cdot f)'}{2f} \le 0
	\quad \mbox{and} \quad
	g'' \le 0
	\qquad \mbox{on } [0,\infty).
  \ee
  Apart from that, we note that with $s_0$ as in (\ref{cinfty_eps}), the first two assumptions in (\ref{struct}) along with
  (\ref{cinfty_eps}) imply that $g'(\ceps) \ge g'(s_0)>0$ in $\Omega\times (0,\tme)$, whence Lemma \ref{lem73}
  combined with (\ref{71.1}) shows that if we take $\rho$ from (\ref{g_Psi_rho}), then
  \bea{74.12}
	\frac{g'(s_0)}{(\cgp)^3} \io \frac{|\nabla \ceps|^4}{\ceps^3}
	&\le& \io \frac{g'(\ceps)}{g^3(\ceps)} |\nabla\ceps|^4 \nn\\
	&\le& (2+\sqrt{3})^2 \io \frac{g(\ceps)}{g'(\ceps)} |D^2 \rho(\ceps)|^2 \nn\\
	&\le& \frac{(2+\sqrt{3})^2}{g'(s_0)} \io g(\ceps)|D^2 \rho(\ceps)|^2
	\qquad \mbox{for all } t\in (0,\tme).
  \eea
  Next, for $i,j\in \{1,2,3\}$ we may again rely on Lemma \ref{lem71} and the concavity and positivity of $g$ and use
  that $(a+b)^2\ge \frac{1}{2}a^2 -b^2$ for $a,b\in\R$ to obtain the pointwise inequality
  \bea{74.122}
	g(\ceps)|\partial_{ij}\rho(\ceps)|^2
	&=& g(\ceps)\Big|\rho'(\ceps) \partial_{ij} \ceps + \rho''(\ceps) \partial_i \ceps \partial_j \ceps \Big|^2 \nn\\
	&\ge& \frac{1}{2} g(\ceps) \rho'^2(\ceps) |\partial_{ij}\ceps|^2
	- g(\ceps) \rho''^2(\ceps) |\partial_i \ceps \partial_j \ceps |^2 \nn\\
	&=& \frac{1}{2g(\ceps)} |\partial_{ij} \ceps|^2
	- \frac{g'^2(\ceps)}{g^3(\ceps)} |\partial_i \ceps \partial_j \ceps |^2 \nn\\
	&\ge& \frac{1}{2\cgp} \cdot \frac{|\partial_{ij}\ceps|^2}{\ceps}
	- \frac{g'^2(0)}{(\cgm)^3} \cdot \frac{|\partial_i \ceps \partial_j \ceps |^2}{\ceps^3}
	\qquad \mbox{in } \Omega\times (0,\tme).
  \eea
  Summarizing, from (\ref{74.12}) and (\ref{74.122}) we thus infer the existence of positive constants $C_1, C_2$ and $C_3$
  such that
  \bas
	\io g(\ceps) |D^2 \rho(\ceps)|^2
	\ge C_1 \io \frac{|\nabla\ceps|^4}{\ceps^3}
	\quad \mbox{and} \quad
	\io g(\ceps) |D^2 \rho(\ceps)|^2
	\ge C_2 \io \frac{|D^2 \ceps|^2}{\ceps}
	-C_3 \io \frac{|\nabla\ceps|^4}{\ceps^3}
  \eas
  for all $t\in (0,\tme)$, which in combination can easily be seen to imply that
  \be{74.123}
	\io g(\ceps) |D^2 \rho(\ceps)|^2
	\ge C_4 \io \frac{|D^2 \ceps|^2}{\ceps} + C_4 \io \frac{|\nabla\ceps|^4}{\ceps^3}
  \ee
  holds for all $t\in (0,\tme)$ if we let $C_4:=\min\{\frac{C_1}{4}, \frac{C_2}{2}, \frac{C_1}{8C_3}\}$.\\
  Since finally $\frac{\partial |\nabla\ceps|^2}{\partial\nu} \le 0$ on $\pO$ according to the convexity of $\Omega$
  (\cite{lions_ARMA}), (\ref{72.1}), (\ref{74.111}) and (\ref{74.123}) imply that 
  \bea{74.2}
	& & \hspace*{-20mm}
	\frac{d}{dt} \bigg\{ \io \neps \ln \neps + \frac{1}{2} \io |\nabla \Psi(\ceps)|^2 \bigg\}
	+ \io \frac{|\nabla\neps|^2}{\neps} 
	+ C_4 \io \frac{|D^2\ceps|^2}{\ceps}
	+ C_4 \io \frac{|\nabla\ceps|^4}{\ceps^3} \nn\\
	&\le& - \frac{1}{2} \io \frac{g'(\ceps)}{g^2(\ceps)} |\nabla \ceps|^2 (\ueps\cdot\nabla \ceps)
	+ \io \frac{1}{g(\ceps)} \Delta\ceps (\ueps\cdot\nabla\ceps)
	\qquad \mbox{for all } t\in (0,\tme).
  \eea
  We now adopt an idea from \cite{chae_kang_lee} and integrate by parts in the rightmost integral herein to see that
  \bas
	\io \frac{1}{g(\ceps)} \Delta\ceps (\ueps\cdot\nabla\ceps)
	&=&
	\io \frac{g'(\ceps)}{g^2(\ceps)} |\nabla\ceps|^2 (\ueps\cdot\nabla\ceps)
	- \io \frac{1}{g(\ceps)} \nabla\ceps \cdot (\nabla\ueps \cdot \nabla\ceps)  \\
	& & - \io \frac{1}{g(\ceps)} \ueps \cdot (D^2\ceps \cdot\nabla\ceps)
  \eas
  for all $t\in (0,\tme)$, where another integration by parts yields
  \bas
	- \io \frac{1}{g(\ceps)} \ueps \cdot (D^2\ceps \cdot\nabla\ceps)
	&=& - \frac{1}{2} \io \frac{1}{g(\ceps)} \ueps \cdot \nabla|\nabla\ceps|^2 \\
	&=& - \frac{1}{2} \io \frac{g'(\ceps)}{g^2(\ceps)} |\nabla\ceps|^2 (\ueps\cdot\nabla\ceps)
	\qquad \mbox{for all } t\in (0,\tme),
  \eas
  because $\nabla\cdot\ueps\equiv 0$.
  Thereby the first integral on the right-hand side of (\ref{74.2}) can be cancelled, so that altogether we obtain
  \bea{74.3}
	\frac{d}{dt} \bigg\{ \io \neps \ln \neps + \frac{1}{2} \io |\nabla \Psi(\ceps)|^2 \bigg\}
	&+& \io \frac{|\nabla\neps|^2}{\neps} 
	+ C_4 \io \frac{|D^2\ceps|^2}{\ceps}
	+ C_4 \io \frac{|\nabla\ceps|^4}{\ceps^3} \nn\\
	&\le& - \io \frac{1}{g(\ceps)} \nabla\ceps \cdot (\nabla\ueps\cdot\nabla\ceps)
  \eea
  for all $t\in (0,\tme)$,
  where by Young's inequality, (\ref{71.1}) and (\ref{cinfty_eps}) we see that with some $C_5>0$ we have
  \bas
	- \io \frac{1}{g(\ceps)} \nabla\ceps \cdot (\nabla\ueps\cdot\nabla\ceps)
	&\le& \frac{C_4}{2} \io \frac{|\nabla\ceps|^4}{\ceps^3}
	+ C_5 \io \frac{\ceps^3}{g^2(\ceps)} |\nabla\ueps|^2 \\
	&\le&  \frac{C_4}{2} \io \frac{|\nabla\ceps|^4}{\ceps^3}
	+ \frac{C_5}{(\cgm)^2} \io \ceps |\nabla\ueps|^2 \\
	&\le&  \frac{C_4}{2} \io \frac{|\nabla\ceps|^4}{\ceps^3}
	+ \frac{C_5 s_0}{(\cgm)^2} \io |\nabla\ueps|^2
	\qquad \mbox{for all } t\in (0,\tme).
  \eas
  The claimed inequalty (\ref{74.1}) thus results from (\ref{74.3}) if we let 
  $\kz:=\max\Big\{ 1 \, , \, \frac{2}{C_4} \, , \, 
  \frac{C_5 s_0}{(\cgm)^2} \Big\}$.
\qed
\subsection{Involving fluid dissipation}\label{sect3.2}		
In order to absorb the source on the right of (\ref{74.1}) appropriately, we recall the following standard energy inequality
for the fluid component (cf.~also \cite{sohr}).
\begin{lem}\label{lem75}
  For each $\eps\in (0,1)$, the solution of (\ref{0eps}) satisfies 
  \be{75.1}
	\frac{1}{2} \frac{d}{dt} \io |\ueps|^2 + \io |\nabla\ueps|^2
	= \io \neps \ueps\cdot \nabla \Phi
	\qquad \mbox{for all } t\in (0,\tme).
  \ee
\end{lem}
\proof
  Testing the third equation in (\ref{0eps}) by $\ueps$ and writing $v_\eps:=\yeps \ueps$, we obtain
  \be{75.2}
	\frac{1}{2} \frac{d}{dt} \io |\ueps|^2 + \io |\nabla\ueps|^2
	= - \io (v_\eps\cdot\nabla)\ueps \cdot \ueps + \io \neps\ueps\cdot\nabla\Phi
	\qquad \mbox{for all } t\in (0,\tme).
  \ee
  Here since $\nabla\cdot\ueps\equiv 0$ and also $\nabla\cdot (I+\eps A)^{-1} \ueps \equiv 0$,
  twice integrating by parts shows that 
  \bas
	\io (v_\eps\cdot\nabla)\ueps \cdot \ueps
	= - \io (\nabla\cdot v_\eps) |\ueps|^2 
	- \frac{1}{2} \io v_\eps \cdot \nabla |\ueps|^2
	= - \frac{1}{2} \io (\nabla\cdot v_\eps) |\ueps|^2 
	=0
  \eas
  for all $t\in (0,\tme)$, whence (\ref{75.2}) implies (\ref{75.1}).
\qed
Now a suitable combination of Lemma \ref{lem74} with Lemma \ref{lem75} yields the following energy-type inequality
which simultaneously involves all the components $\neps,\ceps$ and $\ueps$.
\begin{lem}\label{lem76}
  Let $\Psi$ be as given by (\ref{g_Psi_rho}).
  Then there exist $\kappa>0$ and $K>0$ such that for all $\eps\in (0,1)$, the solution of (\ref{0eps}) satisfies
  \bea{76.1}
	& & \hspace*{-10mm}
	\frac{d}{dt} \bigg\{ \io \neps \ln \neps + \frac{1}{2} \io |\nabla \Psi(\ceps)|^2 + \kappa \io |\ueps|^2 \bigg\}
	+ \frac{1}{K} \bigg\{ \io \frac{|\nabla\neps|^2}{\neps} 
	+ \io \frac{|D^2 \ceps|^2}{\ceps}
	+ \io \frac{|\nabla \ceps|^4}{\ceps^3} + \io |\nabla\ueps|^2 \bigg\} \nn\\[3mm]
	& & \hspace*{68mm}
	\le \ K
	\qquad \mbox{for all } t\in (0,\tme).
  \eea
  In particular, with $\F$ as defined in (\ref{def_energy}) we have
  \bea{776.1}
	-\int_0^\infty \F[\neps,\ceps,\ueps](t) \cdot \phi'(t)dt
	&+& \frac{1}{K} \int_0^\infty \io \bigg\{\frac{|\nabla\neps|^2}{\neps} + \frac{|\nabla \ceps|^4}{\ceps^3}
	+ |\nabla \ueps|^2 \bigg\}(x,t) \cdot \phi(t) dx dt \nn\\[2mm]
	&\le& \F[n_{0\eps}, c_{0\eps}, u_{0\eps}] \cdot \phi(0)
	+ K \int_0^\infty \phi(t) dt
  \eea
  for each nonnegative $\phi\in C_0^\infty([0,\infty))$ and all $\eps\in (0,1)$. 
\end{lem}
\proof
  We first combine (\ref{74.1}) with (\ref{75.1}) to see that with $\kz$ as introduced in Lemma \ref{lem74},
  \bea{76.2}
	\frac{d}{dt} \bigg\{ \io \neps \ln \neps + \frac{1}{2} \io |\nabla \Psi(\ceps)|^2 + \kz \io |\ueps|^2 \bigg\}
	&+& \frac{1}{\kz} \cdot \bigg\{ \io \frac{|\nabla\neps|^2}{\neps} 
	+ \io \frac{|D^2\ceps|^2}{\ceps}
	+ \io \frac{|\nabla \ceps|^4}{\ceps^3} \bigg\} \nn\\[1mm]
	& & + \kz \io |\nabla\ueps|^2 \nn\\[1mm]
	&\le& \kz \io \neps\ueps\cdot\nabla\Phi
	\quad \mbox{for all } t\in (0,\tme).
  \eea
  Here we use the H\"older inequality, (\ref{reg_coeff}) and the continuity of the embedding 
  $W^{1,2}(\Omega)\hra L^6(\Omega)$
  to find $C_1>0$ such that
  \bas
	\kz \io \neps\ueps\cdot\nabla\Phi
	&\le& \kz \|\nabla\Phi\|_{L^\infty(\Omega)} \|\neps\|_{L^\frac{6}{5}(\Omega)} \|\ueps\|_{L^6(\Omega)}  \\
	&\le& C_1\|\neps\|_{L^\frac{6}{5}(\Omega)} \|\nabla\ueps\|_{L^2(\Omega)}
	\qquad \mbox{for all } t\in (0,\tme),
  \eas
  where the Gagliardo-Nirenberg inequality 	
  provides $C_2>0$ and $C_3>0$ such that
  \bas
	\|\neps\|_{L^\frac{6}{5}(\Omega)}
	&=& \|\neps^\frac{1}{2}\|_{L^\frac{12}{5}(\Omega)}^2 \\
	&\le& C_2 \|\nabla\neps^\frac{1}{2}\|_{L^2(\Omega)}^\frac{1}{2} \|\neps^\frac{1}{2}\|_{L^2(\Omega)}^\frac{3}{2}
	+ C_2 \|\neps^\frac{1}{2}\|_{L^2(\Omega)}^2 \\
	&\le& C_3 \cdot \Big\{ \|\nabla\neps^\frac{1}{2}\|_{L^2(\Omega)} + 1 \Big\}^\frac{1}{2}
	\qquad \mbox{for all } t\in (0,\tme),
  \eas
  because $\|\neps^\frac{1}{2}(\cdot,t)\|_{L^2(\Omega)}^2=\io \neps(\cdot,t)=\io n_0$ for all $t\in (0,\tme)$ 
  by (\ref{mass_eps}).
  Twice applying Young's inequality, we hence infer that with some $C_4>0$ and $C_5>0$ we have
  \bas
	\kz \io \neps\ueps\cdot\nabla\Phi
	&\le& C_1 C_3 \cdot \Big\{ \|\nabla\neps^\frac{1}{2}\|_{L^2(\Omega)} + 1 \Big\}^\frac{1}{2}
 	\cdot \|\nabla\ueps\|_{L^2(\Omega)} \\
	&\le& \frac{\kz}{2} \|\nabla\ueps\|_{L^2(\Omega)}^2
	+ C_4 \cdot \Big\{ \|\nabla\neps^\frac{1}{2}\|_{L^2(\Omega)} + 1 \Big\} \\
	&\le& \frac{\kz}{2} \|\nabla\ueps\|_{L^2(\Omega)}^2
	+ \frac{1}{2\kz} \io \frac{|\nabla\neps|^2}{\neps} + C_5
	\qquad \mbox{for all } t\in (0,\tme).
  \eas
  Since $\kz\ge 1$ and hence $\frac{\kz}{2} \ge \frac{1}{2\kz}$, inserted into (\ref{76.2}) this readily yields (\ref{76.1})
  if we let $\kappa:=\kz$ and $K:=\max\{2\kz, C_5\}$.
  Finally, (\ref{776.1}) can be obtained in a straightforward manner 
  on multiplying (\ref{76.1}) by $\phi$ and integrating the resulting
  inequality over $(0,\infty)$, dropping a nonnegative term on its left-hand side.
\qed
In order to derive suitable estimates from this, let us make sure that the energy functional $\F$
therein, when evaluated at the initial time, approaches its expected limit as $\eps\searrow 0$.
Since in this respect the integrals involving $n_{0\eps}$ and $u_{0\eps}$ clearly have the desired behavior 
due to (\ref{I1}) and (\ref{I3}), this actually reduces to proving the following lemma, which for later purpose
asserts a slightly more general statement.
\begin{lem}\label{lem70}
  Let $\Psi$ be as in (\ref{g_Psi_rho}), and suppose that $(\varphi_j)_{j\in\N} \subset C^2(\bar\Omega;(0,\infty))$
  and $\varphi:\Omega\to [0,\infty)$ are such that $\sqrt{\varphi}\in W^{1,2}(\Omega)$ and
  \be{70.1}
	\sqrt{\varphi_j} \to \sqrt{\varphi}
	\quad \mbox{in } W^{1,2}(\Omega)
	\mbox{ and a.e.~in } \Omega
	\quad \mbox{as } j\to\infty
  \ee
  as well as $\|\varphi_j\|_{L^\infty(\Omega)} \le s_0$ with $s_0$ given by (\ref{cinfty_eps}).
  Then $\Psi(\varphi)\in W^{1,2}(\Omega)$ and
  \be{70.2}
	\Psi(\varphi_j) \to \Psi(\varphi)
	\quad \mbox{in } W^{1,2}(\Omega)
	\qquad \mbox{as } j\to\infty.
  \ee
\end{lem}
\proof
  Since (\ref{71.1}) and the inequality $\varphi_j\le s_0$ ensure that 
  \bas
	\varphi_j \Psi'^2(\varphi_j) = \frac{\varphi_j}{g(\varphi_j)} \le \frac{1}{\cgm}
	\quad \mbox{in } \Omega
	\qquad \mbox{for all } \eps\in (0,1),
  \eas
  our assumption that $\varphi_j\to\varphi $ a.e.~in $\Omega$ as $j\to\infty$
  entails that
  $\varphi_j \Psi'^2(\varphi_j) \wsto \varphi \Psi'^2(\varphi)$ in $L^\infty(\Omega)$ as $j\to\infty$.
  Combined with the fact that by (\ref{70.1}) we have
  $|\nabla\sqrt{\varphi_j}|^2 \to |\nabla \sqrt{\varphi}|^2$
  in $L^1(\Omega)$ as $j\to\infty$, this shows that
  \bas
	\io |\nabla \Psi(\varphi_j)|^2
	= \io \varphi_j \Psi'^2(\varphi_j) |\nabla \sqrt{\varphi_j}|^2
	\to \io \varphi \Psi'^2(\varphi) |\nabla\sqrt{\varphi}|^2
	= \io |\nabla \Psi(\varphi)|^2
  \eas
  as $j\to\infty$.
  Since the estimates in Lemma \ref{lem70} along with (\ref{70.1}) and the dominated convergence theorem 
  readily imply that $\Psi(\varphi_j)\to\Psi(\varphi)$ in $L^1(\Omega)$ as $j\to\infty$, this proves (\ref{71.1}).
\qed
We can thereby draw the following consequence of Lemma \ref{lem76}.
\begin{lem}\label{lem77}
  There exists $C>0$ such that with $\Psi$ as in (\ref{g_Psi_rho}) we have
  \be{77.01}
	\io \neps(\cdot,t)\ln \neps(\cdot,t) + 
	\io |\nabla \Psi(\ceps(\cdot,t))|^2 + 
	\io |\ueps(\cdot,t)|^2 \le C
	\qquad \mbox{for all } t\in (0,\tme)
  \ee
  and
  \be{77.1}
	\int_0^T \io \frac{|\nabla\neps|^2}{\neps}
	+ \int_0^T \io \frac{|D^2 \ceps|^2}{\ceps}
	+ \int_0^T \io \frac{|\nabla \ceps|^4}{\ceps^3} 
	+ \int_0^T \io |\nabla\ueps|^2
	\le C\cdot (T+1)
	\qquad \mbox{for all } T\in (0,\tme)
  \ee
  whenever $\eps\in (0,1)$.
\end{lem}
\proof 
  With $\kappa>0$ and $K>0$ as provided by Lemma \ref{lem76}, an application of the latter shows that if we take
  $\Psi$ from (\ref{g_Psi_rho}), then for each 
  $\eps\in (0,1)$,
  \bas
	y_\eps(t):= \F[\neps,\ceps,\ueps](t) \equiv
	\io \bigg\{ \neps\ln \neps + \frac{1}{2}|\nabla\Psi(\ceps)|^2 + \kappa |\ueps|^2 \bigg\} (\cdot,t),
	\qquad t\in [0,\tme),
  \eas
  satisfies
  \be{77.2}
	y_\eps'(t) + \frac{1}{K} h_\eps(t) \le K
	\qquad \mbox{for all } t\in (0,\tme),
  \ee
  where 
  \bas
	h_\eps(t):=\io \bigg\{ \frac{|\nabla\neps|^2}{\neps} + 	
	\frac{|D^2 \ceps|^2}{\ceps} 
	+ \frac{|\nabla\ceps|^4}{\ceps^3} + |\nabla\ueps|^2 \bigg\}
	(\cdot,t)
	\qquad \mbox{for } t\in (0,\tme).
  \eas
  In order to estimate $y_\eps(t)$ in terms of $h_\eps(t)$, we first use the Poincar\'e inequality to find $C_1>0$ such that
  \be{77.3}
	\io |\ueps|^2 \le C_1 \io |\nabla\ueps|^2
	\qquad \mbox{for all } t\in (0,\tme),
  \ee
  and next recall the definitions of $\Psi$ and $g$ in (\ref{g_Psi_rho}) as well as (\ref{71.1}) and (\ref{cinfty_eps}) 
  to see employing Young's inequality that
  \bea{77.4}
	\frac{1}{2} \io |\nabla\Psi(\ceps)|^2
	&=& \frac{1}{2} \io \frac{|\nabla\ceps|^2}{g(\ceps)} \nn\\
	&\le& \io \frac{|\nabla\ceps|^4}{\ceps^3}
	+ \frac{1}{16} \io \frac{\ceps^3}{g^2(\ceps)} \nn\\
	&\le& \io \frac{|\nabla\ceps|^4}{\ceps^3}
	+ \frac{1}{16(\cgm)^2} \io \ceps \nn\\
	&\le& \io \frac{|\nabla\ceps|^4}{\ceps^3}
	+ \frac{s_0|\Omega|}{16(\cgm)^2}
	\qquad \mbox{for all } t\in (0,\tme).
  \eea
  We finally make use of the elementary inequality $z\ln z \le \frac{3}{2} z^\frac{5}{3}$, valid for all $z\ge 0$,
  and invoke the Gagliardo-Nirenberg inequality together with (\ref{mass_eps}) to infer that with some 
  $C_2>0$ and $C_3>0$ we have
  \bas
	\io \neps\ln \neps
	&\le& \frac{3}{2} \io \neps^\frac{5}{3} \\
	&=& \frac{3}{2} \|\neps^\frac{1}{2}\|_{L^\frac{10}{3}(\Omega)}^\frac{10}{3}  \\
	&\le& C_2 \|\nabla\neps^\frac{1}{2}\|_{L^2(\Omega)}^2 \|\neps^\frac{1}{2}\|_{L^2(\Omega)}^\frac{4}{3}
	+ C_2 \|\neps^\frac{1}{2}\|_{L^2(\Omega)}^\frac{10}{3} \\
	&\le& C_3 \io \frac{|\nabla\neps|^2}{\neps} + C_3
  \eas
  for all $t\in (0,\tme)$.
  In conjunction with (\ref{77.3}) and (\ref{77.4}), this provides $C_4>0$ such that
  \bas
	y_\eps(t) \le C_4 h_\eps(t) + C_4
	\qquad \mbox{for all } t\in (0,\tme),
  \eas
  so that (\ref{77.2}) implies that $y_\eps$ satisfies the ODI
  \be{77.5}
	y_\eps'(t) + \frac{1}{2K} h_\eps(s) + \frac{1}{2KC_4} y_\eps(t)
	\le C_5:=K+\frac{1}{2K}
	\qquad \mbox{for all } t\in (0,\tme).
  \ee
  This firstly warrants that 
  \be{77.6}
	y_\eps(t) \le C_6:=\max \Big\{ \sup_{\eps\in (0,1)} y_\eps(0) \, , \, 2KC_4 C_5\Big\}
	\qquad \mbox{for all } t\in (0,\tme)
  \ee
  and thus proves (\ref{77.01}), because $\sup_{\eps\in (0,1)} y_\eps(0)$ is finite thanks to (\ref{I1}), (\ref{I2}),
  (\ref{I3}) and Lemma \ref{lem70}.
  Secondly, another integration of (\ref{77.5}) thereupon shows that
  \bas
	\frac{1}{2K} \int_0^T h_\eps(t)dt \le y_\eps(0) + C_5 T \le C_6 + C_5 T
	\qquad \mbox{for all } T\in (0,\tme),
  \eas
  which in view of the definition of $h_\eps$ establishes (\ref{77.1}).
\qed
\subsection{Global existence in the regularized problems}\label{sect3.3}
In light of the fact that our specific choice (\ref{Feps})
of $F_\eps$ warrants that $[0,\infty) \ni s \mapsto sF_\eps'(s)$ is bounded
for any fixed $\eps\in (0,1)$,
the bound for $\nabla \ceps$ in $L^4_{loc}(\bar\Omega\times [0,\tme))$ implied by Lemma \ref{lem77} and (\ref{cinfty_eps})
is sufficient to guarantee that each of our approximate solutions is indeed global in time:
\begin{lem}\label{lem777}
  For all $\eps\in (0,1)$, the solution of (\ref{0eps}) is global in time; that is, we have $\tme=\infty$.
\end{lem}
\proof
  Assuming that $\tme$ be finite for some $\eps\in (0,1)$, we first note that as a particular
  consequence of Lemma \ref{lem77} and (\ref{cinfty_eps}) we can then find $C_1>0$ and $C_2>0$ such that
  \be{777.01}
	\int_0^{\tme} \io |\nabla \ceps|^4 \le C_1
	\qquad \mbox{and} \qquad
	\io |\ueps(\cdot,t)|^2 \le C_2
	\quad \mbox{for all } t\in (0,\tme).
  \ee
  To deduce a contradiction from this, we 
  firstly multiply the equation for $\neps$ in (\ref{0eps}) by $\neps^3$, integrate by parts and use Young's inequality
  together with the fact that $\neps F_\eps'(\neps) \le \frac{1}{\eps}$ by (\ref{Feps1}) 
  to obtain $C_3>0$, as all constants below possibly depending on $\eps$, such that
  \bas
	\frac{1}{4} \frac{d}{dt} \io \neps^4 + 3 \io \neps^2 |\nabla\neps|^2 	
	\le \io \neps^2 |\nabla \neps|^2 + \io |\nabla \ceps|^4 + C_3 \io \neps^4 
	\qquad \mbox{for all } t\in (0,\tme),
  \eas
  so that thanks to the first inequality in (\ref{777.01}) we see that
  \be{777.1}
	\io \neps^4(\cdot,t)\le C_4
	\qquad \mbox{for all } t\in (0,\tme)
  \ee
  with some $C_4>0$. \\
  We next observe that $D(1+\eps A)=W^{2,2}(\Omega)\cap W_{0,\sigma}^{1,2}(\Omega) \hra L^\infty(\Omega)$,
  according to the second estimate in (\ref{777.01}) there exist $C_5>0$ and $C_6>0$ such that
  $\veps:=\yeps \ueps$ satisfies
  \be{777.33}
	\|\veps(\cdot,t)\|_{L^\infty(\Omega)} 
	= \big\| (1+\eps A)^{-1} \ueps(\cdot,t) \big\|_{L^2(\Omega)}
	\le C_4 \|\ueps(\cdot,t)\|_{L^2(\Omega)} \le C_6
	\qquad \mbox{for all } t\in (0,\tme).
  \ee
  Therefore, testing the projected Stokes equation 
  $u_{\eps t}+A \ueps = h_\eps(x,t):= \proj[ - (\veps\cdot\nabla)\ueps + \neps \nabla \Phi]$ by $A\ueps$ shows that
  \bas
	\frac{1}{2} \frac{d}{dt} \io |A^\frac{1}{2} \ueps|^2 + \io |A\ueps|^2
	&=& \io A\ueps \cdot h_\eps \\
	&\le& \io |A\ueps|^2 + \frac{1}{4} \io |h_\eps|^2 \\
	&\le& \io |A\ueps|^2 
	+ \frac{1}{2} \io |(\veps\cdot\nabla)\ueps|^2
	+ \frac{1}{2} \io |\neps\nabla\Phi|^2 \\
	&\le& \io |A\ueps|^2 
	+ C_7 \cdot \Big\{ \io |\nabla\ueps|^2 + \io \neps^2 \Big\}
	\qquad \mbox{for all } t\in (0,\tme)
  \eas	
  with some $C_7>0$, because $\|\proj\varphi\|_{L^2(\Omega)} \le \|\varphi\|_{L^2(\Omega)}$ for all
  $\varphi\in L^2(\Omega)$.
  As $\io |A^\frac{1}{2}\varphi|^2=\io |\nabla\varphi|^2$ for all $\varphi\in D(A)$, in view of (\ref{777.1})
  this implies the existence of $C_8>0$ fulfilling
  \be{777.3}
	\io |\nabla\ueps(\cdot,t)|^2 \le C_8
	\qquad \mbox{for all } t\in (0,\tme).
  \ee
  Along with (\ref{777.33}) and again (\ref{777.1}), this in turn provides $C_9>0$ such that 
  $\|h_\eps(\cdot,t)\|_{L^2(\Omega)} \le C_9$ for all $t\in (0,\tme)$.
  Thus if we pick an arbitrary $\alpha\in (\frac{3}{4},1)$, then known smoothing properties of the Stokes semigroup
  (\cite[p. 201]{giga1986}) entail that for some $C_{10}>0$ we have
  \bas
	\|A^\alpha \ueps(\cdot,t)\|_{L^2(\Omega)}
	&=& \bigg\| A^\alpha e^{-tA} u_{0\eps} + \int_0^t A^\alpha e^{-(t-s)A} h_\eps(\cdot,s) ds \bigg\|_{L^2(\Omega)} \\
	&\le& C_{10} t^{-\alpha} \|u_{0\eps}\|_{L^2(\Omega)}
	+ C_{10} \int_0^t (t-s)^{-\alpha} \|h_\eps(\cdot,s)\|_{L^2(\Omega)} ds \\
	&\le& C_{10} t^{-\alpha} \|u_{0\eps}\|_{L^2(\Omega)}
	+ \frac{C_9 C_{10} \tme^{1-\alpha}}{1-\alpha}
	\qquad \mbox{for all } t\in (0,\tme).
  \eas
  Since also $D(A^\alpha)$ is continuously embedded into $L^\infty(\Omega)$ due to our choice of $\alpha$
  (\cite[Thm.~1.6.1]{henry}, \cite{giga1981_the_other}), we thereby obtain
  $C_{11}>0$ and $C_{12}>0$ such that with $\tau:=\frac{1}{2} \tme$,
  \be{777.5}
	\|\ueps(\cdot,t)\|_{L^\infty(\Omega)}
	\le C_{11} \|A^\alpha \ueps(\cdot,t)\|_{L^2(\Omega)} \le C_{12}
	\qquad \mbox{for all } t\in (\tau,\tme).
  \ee
  Thereafter, standard smoothing estimates for the Neumann heat semigroup (\cite{quittner_souplet}, \cite{win_JDE}),
  an application of the H\"older inequality
  and (\ref{cinfty_eps}), (\ref{Feps1}),
  (\ref{777.1}), (\ref{777.5}) and (\ref{777.01}) yield positive constants $C_{13}, C_{14}$ and $C_{15}$ satisfying
  \bea{777.7}
	\|\nabla\ceps(\cdot,t)\|_{L^4(\Omega)}
	&=& \bigg\| \nabla e^{(t-\tau)\Delta} \ceps(\cdot,\tau)
	- \int_\tau^t \nabla e^{(t-s)\Delta} \Big\{ F_\eps(\neps) f(\ceps) + \ueps\cdot\nabla\ceps \Big\} (\cdot,s) ds 
	\bigg\|_{L^4(\Omega)} \nn\\
	&\le& C_{13} (t-\tau)^{-\frac{1}{2}} \|\ceps(\cdot,\tau)\|_{L^4(\Omega)} \nn\\
	& & + C_{13} \int_\tau^t (t-s)^{-\frac{1}{2}} 
	\bigg\{ \Big\| F_\eps(\neps(\cdot,s)) f(\ceps(\cdot,s))\Big\|_{L^4(\Omega)}
	+ \Big\|\ueps(\cdot,s)\cdot\nabla\ceps(\cdot,s)\Big\|_{L^4(\Omega)} \bigg\} ds \nn\\
	&\le& C_{13} (t-\tau)^{-\frac{1}{2}} \|\ceps(\cdot,\tau)\|_{L^4(\Omega)} \nn\\
	& & + C_{14} \int_\tau^t (t-s)^{-\frac{1}{2}} ds
	+ C_{14} \int_\tau^t (t-s)^{-\frac{1}{2}} \|\nabla \ceps(\cdot,s)\|_{L^4(\Omega)} ds \nn\\
	&\le& C_{13} (t-\tau)^{-\frac{1}{2}} \|\ceps(\cdot,\tau)\|_{L^4(\Omega)} 
	+ 2C_{14} \tme^\frac{1}{2} \\
	& & + C_{14} \bigg\{ \int_0^{\tme} \sigma^{-\frac{2}{3}} d\sigma \bigg\}^\frac{3}{4}
	\cdot \bigg\{ \int_0^{\tme} \|\nabla\ceps(\cdot,s)\|_{L^4(\Omega)}^4 ds \bigg\}^\frac{1}{4} \nn\\[2mm]
	&\le& C_{15}
	\qquad \mbox{for all } t\in (2\tau,\tme).
  \eea
  Similarly, combining (\ref{cinfty_eps}), (\ref{Feps1}), (\ref{777.1}) and (\ref{777.5}) shows that
  there exist $C_{16}>0, C_{17}>0$ and $C_{18}>0$ such that
  \bas
	\|\neps(\cdot,t)\|_{L^\infty(\Omega)}
	&=& \bigg\| e^{(t-\tau)\Delta} \neps(\cdot,\tau) 
	- \int_\tau^t e^{(t-s)\Delta} \nabla \cdot \Big\{ F_\eps(\neps(\cdot,s)) f(\ceps(\cdot,s)) 
	+ \neps(\cdot,s)\ueps(\cdot,s) \Big\} ds \bigg\|_{L^\infty(\Omega)} \\
	&\le& C_{16} (t-\tau)^{-\frac{3}{2}} \|\neps(\cdot,\tau)\|_{L^1(\Omega)} \\
	& & + C_{16} \int_\tau^t (t-s)^{-\frac{1}{2} -\frac{3}{2} \cdot \frac{1}{4}} 
	\|\neps(\cdot,s)\|_{L^4(\Omega)} \Big\{ 1 + \|\ueps(\cdot,s)\|_{L^\infty(\Omega)} \Big\} ds \\[2mm]
	&\le& C_{17}
	\qquad \mbox{for all } t\in (2\tau,\tme).
  \eas
  Together with (\ref{777.7}) and (\ref{777.5}), this contradicts the extensibility criterion (\ref{ext_crit})
  in Lemma \ref{lem90} and thereby entails that actually $\tme=\infty$, as claimed.
\qed
\subsection{Further a priori estimates. Time regularity}\label{sect3.4}
By interpolation, the estimates from Lemma \ref{lem77} imply bounds for further spatio-temporal integrals.
\begin{lem}\label{lem80}
  There exists $C>0$ such that for each $\eps\in (0,1)$ we have
  \be{80.1}
	\int_0^T \io \neps^\frac{5}{3} 
	\le C\cdot (T+1) \qquad \mbox{for all } T>0
  \ee
  and
  \be{80.2}
	\int_0^T \io |\nabla \neps|^\frac{5}{4}
	\le C\cdot (T+1) \qquad \mbox{for all } T>0
  \ee
  as well as
  \be{80.3}
	\int_0^T \io |\ueps|^\frac{10}{3}
	\le C\cdot (T+1) \qquad \mbox{for all } T>0.
  \ee
\end{lem}
\proof
  According to Lemma \ref{lem77}, there exists $C_1>0$ such that
  \be{80.4}
	\int_0^T \io \frac{|\nabla \neps|^2}{\neps} \le C_1\cdot(T+1)
	\qquad \mbox{for all } T>0.
  \ee
  Thus, invoking the Gagliardo-Nirenberg inequality along with (\ref{mass_eps}) 
  we obtain $C_2>0$ and $C_3>0$ such that
  \bea{80.5}
	\int_0^T \io \neps^\frac{5}{3}
	&=& \int_0^T \|\neps^\frac{1}{2}(\cdot,t)\|_{L^\frac{10}{3}(\Omega)}^\frac{10}{3} dt \nn\\
	&\le& C_2 \int_0^T \Big\{ \|\nabla \neps^\frac{1}{2}(\cdot,t)\|_{L^2(\Omega)}^2
	\cdot \|\neps^\frac{1}{2}(\cdot,t)\|_{L^2(\Omega)}^\frac{4}{3}
	+ \|\neps^\frac{1}{2}(\cdot,t)\|_{L^2(\Omega)}^\frac{10}{3} \Big\} dt \nn\\
	&\le& C_2 \cdot \Big\{ \frac{C_1\cdot(T+1)}{4} \cdot \|n_0\|_{L^1(\Omega)}^\frac{2}{3} 
	+ \|n_0\|_{L^1(\Omega)}^\frac{5}{3} T \Big\} \nn\\[2mm]
	&\le& C_3\cdot (T+1)
	\qquad \mbox{for all } T>0.
  \eea
  Employing the H\"older inequality, again by (\ref{mass_eps}) we furthermore conclude from (\ref{80.4})
  together with (\ref{80.5}) that
  \bas		
	\int_0^T \io |\nabla\neps|^\frac{5}{4}
	&=& \int_0^T \io \frac{|\nabla\neps|^\frac{5}{4}}{\neps^\frac{5}{8}} \cdot \neps^\frac{5}{8} \nn\\
	&\le& \bigg( \int_0^T \io \frac{|\nabla\neps|^2}{\neps} \bigg)^\frac{5}{8} \cdot
	\bigg(\int_0^T \io \neps^\frac{5}{3}\bigg)^\frac{3}{8} \nn\\
	&\le& C_1^\frac{5}{8} C_3^\frac{3}{8} \cdot (T+1)
	\qquad \mbox{for all } T>0.
  \eas
  As once more relying on Lemma \ref{lem77} we can find $C_4>0$ and $C_5>0$ such that
  \bas
	\io |\ueps(\cdot,t)|^2 \le C_4
	\quad \mbox{for all } t>0
	\qquad \mbox{and} \qquad
	\int_0^T \io |\nabla\ueps|^2 \le C_5 \cdot (T+1)
	\qquad \mbox{for all } T>0,
  \eas
  upon another application of the Gagliardo-Nirenberg inequality in precisely the same way as in (\ref{80.5})
  we see that with some $C_6>0$ we have
  \bas
	\int_0^T \|\ueps(\cdot,t)\|_{L^\frac{10}{3}(\Omega)}^\frac{10}{3} dt
	&\le& C_6 \int_0^T \Big\{ \|\nabla \ueps(\cdot,t)\|_{L^2(\Omega)}^2
	\cdot \|\ueps(\cdot,t)\|_{L^2(\Omega)}^\frac{4}{3}
	+ \|\ueps(\cdot,t)\|_{L^2(\Omega)}^\frac{10}{3} \Big\} dt \nn\\
	&\le& C_4^\frac{4}{3} C_5 C_6 \cdot (T+1)
	+ C_4^\frac{10}{3} C_6 T
	\qquad \mbox{for all } T>0,
  \eas
  whereby the proof is completed.
\qed
In a straightforward manner, from Lemma \ref{lem77} and Lemma \ref{lem80} we can moreover deduce
certain regularity features of the time derivatives in (\ref{0eps}).
Since in Lemma \ref{lem82}
these will mainly be used to warrant pointwise convergence 
we refrain from pursuing here the question which are the smallest spaces within which such derivative bounds can be obtained.
\begin{lem}\label{lem81}
  There exists $C>0$ such that
  \be{81.1}
	\int_0^T \|n_{\eps t}(\cdot,t)\|_{(W^{1,10}(\Omega))^\star}^\frac{10}{9} dt
	\le C\cdot (T+1) \qquad \mbox{for all } T>0
  \ee
  and
  \be{81.2}
	\int_0^T \|(\sqrt{c_{\eps})_t}(\cdot,t)\|_{(W^{1,\frac{5}{2}}(\Omega))^\star}^\frac{5}{3} dt
	\le C\cdot (T+1) \qquad \mbox{for all } T>0
  \ee
  as well as
  \be{81.3}
	\int_0^T \|u_{\eps t}(\cdot,t)\|_{(W_{0,\sigma}^{1,5}(\Omega))^\star}^\frac{5}{4} dt
	\le C\cdot (T+1) \qquad \mbox{for all } T>0.
  \ee
\end{lem}
\proof
  For arbitrary $t>0$ and $\varphi\in C^\infty(\bar\Omega)$, multiplying the first equation in (\ref{0eps}) by $\varphi$,
  integrating by parts and
  using the H\"older inequality we obtain
  \bas
	\bigg| \io n_{\eps t}(\cdot,t)\varphi \bigg|
	&=& \bigg| - \io \nabla\neps\cdot \nabla\varphi 
	+ \io \neps F_\eps'(\neps) \chi(\ceps) \nabla\ceps\cdot\nabla \varphi
	+\io \neps\ueps \cdot\nabla\varphi \bigg| \\
	&\le& \bigg\{ \|\nabla\neps\|_{L^\frac{10}{9}(\Omega)} 
	+ \|\neps F_\eps'(\neps) \chi(\ceps) \nabla\ceps\|_{L^\frac{10}{9}(\Omega)}
	+ \|\neps \ueps\|_{L^\frac{10}{9}(\Omega)} \bigg\}
	\cdot \|\varphi\|_{W^{1,10}(\Omega)},
  \eas
  so that with some $C_1>0$ we have
  \bea{81.4}
	\int_0^T \|n_{\eps t}(\cdot,t)\|_{(W^{1,10}(\Omega))^\star}^\frac{10}{9} dt
	&\le& \int_0^T \bigg\{ \|\nabla\neps\|_{L^\frac{10}{9}(\Omega)} 
	+ \|\neps F_\eps'(\neps) \chi(\ceps) \nabla\ceps\|_{L^\frac{10}{9}(\Omega)}
	+ \|\neps \ueps\|_{L^\frac{10}{9}(\Omega)} \bigg\}^\frac{10}{9} dt \nn\\
	&\le& C_1\int_0^T \io |\nabla \neps|^\frac{10}{9}
	+ C_1 \int_0^T \io |\neps \nabla \ceps |^\frac{10}{9}
	+ C_1 \int_0^T \io |\neps\ueps|^\frac{10}{9}
  \eea
  for all $T>0$,
  because $F_\eps'(\neps)\le 1$ by (\ref{Feps1}) and $\chi(\ceps)\le \|\chi\|_{L^\infty((0,M))}$ with
  $M:=\|c_0\|_{L^\infty(\Omega)}$ according to (\ref{cinfty_eps}).
  Here several applications of Young's inequality show that
  \bas
	\int_0^T \io |\nabla\neps|^\frac{10}{9} \le \int_0^T \io |\nabla\neps|^\frac{5}{4} + |\Omega| T
	\qquad \mbox{for all } T>0
  \eas
  and
  \bas
	\int_0^T \io |\neps\nabla \ceps|^\frac{10}{9}
	\le \int_0^T \io \neps^\frac{5}{3}
	+ \int_0^T \io |\nabla \ceps|^\frac{10}{3}
	\le \int_0^T \io \neps^\frac{5}{3}
	+ \int_0^T \io |\nabla \ceps|^4 + |\Omega|T
	\qquad \mbox{for all } T>0
  \eas
  as well as
  \bas
	\int_0^T \io |\neps\ueps|^\frac{10}{9}
	\le \int_0^T \io \neps^\frac{5}{3} + \int_0^T \io |\ueps|^\frac{10}{3}
	\qquad \mbox{for all } T>0,
  \eas
  whence in light of Lemma \ref{lem80}, Lemma \ref{lem77} and (\ref{cinfty_eps}), (\ref{81.1}) results from (\ref{81.4}).\abs
  Likewise, given any $\varphi\in C^\infty(\bar\Omega)$ we may test the second equation in (\ref{0eps}) against 
  $\frac{\varphi}{\sqrt{\ceps}}$ to see that 
  \bas
	\bigg| \io (\sqrt{\ceps})_t(\cdot,t) \varphi \bigg|
	&=& \bigg| - \frac{1}{2} \io \frac{\nabla\ceps}{\sqrt{\ceps}}\cdot\nabla\varphi 
	+ \frac{1}{4} \io \frac{|\nabla\ceps|^2}{\sqrt{\ceps}^3} \varphi
	- \frac{1}{2} \io F_\eps(\neps) \frac{f(\ceps)}{\sqrt{\ceps}} \varphi
	+ \io \sqrt{\ceps} \ueps \cdot \nabla\varphi \bigg| \\
	& & \hspace*{-20mm}
	\le \ \bigg\{ \frac{1}{2} \Big\| \frac{\nabla\ceps}{\sqrt{\ceps}} \Big\|_{L^\frac{5}{3}(\Omega)}
	+ \frac{1}{4} \Big\| \frac{|\nabla\ceps|^2}{\sqrt{\ceps}^3} \Big\|_{L^\frac{5}{3}(\Omega)}
	+ \Big\|F_\eps(\neps)\frac{f(\ceps)}{\sqrt{\ceps}}\Big\|_{L^\frac{5}{3}(\Omega)}
	+ \|\sqrt{\ceps}\ueps\|_{L^\frac{5}{3}(\Omega)} \bigg\} \cdot \|\varphi\|_{W^{1,\frac{5}{2}}(\Omega)}
  \eas
  for all $t>0$, and that hence by Young's inequality we can find $C_2>0$ such that
  \bas
	\int_0^T \|(\sqrt{\ceps})_t(\cdot,t)\|_{(W^{1,\frac{5}{2}}(\Omega))^\star}^\frac{5}{3} dt
	&\le& C_2 \int_0^T \io \ceps^{-\frac{5}{6}} |\nabla\ceps|^\frac{5}{3}
	+ C_2 \int_0^T \io \ceps^{-\frac{5}{2}} |\nabla\ceps|^\frac{10}{3} \nn\\
	& & + C_2 \int_0^T \io \neps^\frac{5}{3}
	+ C_2 \int_0^t \io |\ueps|^\frac{5}{3} \\
	&\le& C_2 \int_0^T \io \frac{|\nabla\ceps|^4}{\ceps^3} 
	+ C_2 M^\frac{5}{7} |\Omega|T \nn\\
	& & + C_2 \int_0^T \io \frac{|\nabla\ceps|^4}{\ceps^3} 
	+ C_2 |\Omega|T \nn\\
	& & + C_2 \int_0^T \io \neps^\frac{5}{3} + C_2\int_0^T \io |\ueps|^\frac{10}{3}
	+ C_2 |\Omega|T
  \eas
  for all $T>0$, since $F_\eps(\neps)\le \neps$ due to (\ref{Feps1}), and since once more by (\ref{cinfty_eps})
  we have $\ceps\le s_0$ and $\frac{f(\ceps)}{\sqrt{\ceps}}\le \|f\|_{C^\frac{1}{2}([0,s_0])}$ in $\Omega\times (0,\infty)$.
  Again by Lemma \ref{lem80} and Lemma \ref{lem77}, this imples (\ref{81.2}).\\
  Finally, given $\varphi\in C_{0,\sigma}^\infty(\Omega;\R^3)$ we infer from the third equation in (\ref{0eps}) that
  \bea{81.999}
	\bigg| \io u_{\eps t}(\cdot,t) \cdot \varphi \bigg|
	&=& \bigg| - \io \nabla\ueps \cdot \nabla \varphi
	- \io (\yeps \ueps \mult \ueps ) \cdot \nabla\varphi
	+ \io \neps \nabla\Phi \cdot\nabla\varphi \bigg| \nn\\
	&\le& \bigg\{ \|\nabla\ueps\|_{L^\frac{5}{4}(\Omega)}
	+ \|\yeps \ueps \mult \ueps \|_{L^\frac{5}{4}(\Omega)}
	+ \|\neps\nabla\Phi\|_{L^\frac{5}{4}(\Omega)} \bigg\} \cdot \|\varphi\|_{W^{1,5}(\Omega)}
  \eea
  for all $t>0$, 
  where for $v\in\R^3$ and $w\in\R^3$ we have defined the matrix $v\mult w$ by letting its components be given by 
  $(v\mult w)_{ij}:= v_i w_j$ for $i,j\in\{1,2,3\}$.
  In view of Young's inequality, (\ref{81.999}) implies that there exists $C_3>0$ fulfilling
  \bas
	\int_0^T \|u_{\eps t}(\cdot,t)\|_{(W_{0,\sigma}^{1,5}(\Omega))^\star}^\frac{5}{4} dt
	&\le& C_3 \int_0^T \io |\nabla \ueps|^\frac{5}{4}
	+ C_3 \int_0^T \io |\yeps \ueps \mult \ueps|^\frac{5}{4}
	+ C_3 \int_0^T \io \neps^\frac{5}{4} \\
	&\le& C_3 \int_0^T \io |\nabla\ueps|^2
	+ C_3 \int_0^T \io |\yeps \ueps|^2
	+ C_3 \int_0^T \io |\ueps|^\frac{10}{3} \\
	& & + C_3 \int_0^T \io \neps^\frac{5}{3}
	+ 2C_3|\Omega|T
	\qquad \mbox{for all } T>0,
  \eas
  because $\nabla\Phi\in L^\infty(\Omega)$.
  Since $\|\yeps v\|_{L^2(\Omega)} \le \|v\|_{L^2(\Omega)}$ for all $v\in L^2_\sigma(\Omega)$ and hence 
  $\int_0^T \io |\yeps \ueps|^2 \le \int_0^T \io |\ueps|^\frac{10}{3} + |\Omega|T$ for all $T>0$, (\ref{81.3})
  results from this upon another application of Lemma \ref{lem80} and Lemma \ref{lem77}.
\qed
\mysection{Passing to the limit. Proof of Theorem \ref{theo_global}}\label{sect4}
With the above compactness properties at hand, by means of a standard extraction procedure we can now derive the
following lemma which actually contains our main existence result already.
\begin{lem}\label{lem82}
  There exists $(\eps_j)_{j\in\N} \subset (0,1)$ 	
  such that $\eps_j\searrow 0$ as $j\to\infty$, and such that as $\eps=\eps_j\searrow 0$ we have
  \begin{eqnarray}
	\neps\to n
	& & 
	\mbox{in } L^\frac{5}{3}_{loc}(\bar\Omega\times [0,\infty))
	\ \mbox{and a.e.~in } \Omega\times (0,\infty), \label{82.1} \\
	\nabla\neps\wto \nabla n
	& & 
	\mbox{in } L^\frac{5}{4}_{loc}(\bar\Omega\times [0,\infty)), \label{82.2} \\
	\ceps\to c
	& & 
	\mbox{a.e.~in } \Omega\times (0,\infty), \label{82.3} \\
	\ceps\wsto c
	& & 
	\mbox{in } L^\infty(\Omega\times (0,\infty)), \label{82.4} \\
	\nabla\ceps^\frac{1}{4} \wto \nabla c^\frac{1}{4}
	& & 
	\mbox{in } L^4_{loc}(\bar\Omega\times [0,\infty)), \label{82.107} \\
	\ueps\to u
	& & 
	\mbox{in } L^2_{loc}(\bar\Omega\times [0,\infty))
	\ \mbox{and a.e.~in } \Omega\times (0,\infty), \label{82.6} \\
	\ueps\wsto u
	& & 
	\mbox{in } L^\infty([0,\infty);L^2_\sigma(\Omega)), \label{82.666} \\
	\ueps \wto u
	& & 
	\mbox{in } L^\frac{10}{3}_{loc}(\bar\Omega\times [0,\infty)) \qquad \mbox{and} \label{82.66} \\
	\nabla\ueps \wto \nabla u
	& &
	\mbox{in } L^2_{loc}(\bar\Omega\times [0,\infty)) \label{82.7}
  \end{eqnarray}
  with some limit functions $n,c$ and $u$ such that $(n,c,u)$ is a global weak solution of (\ref{0}), (\ref{0i}), (\ref{0b})
  in the sense of Definition \ref{defi_weak}.
\end{lem}
\proof 
  Using the pointwise identity
  \bas
	\partial_{ij} \sqrt{\ceps} = \frac{1}{2\sqrt{\ceps}} \partial_{ij} \ceps
	- \frac{1}{4\sqrt{\ceps}^3} \partial_j \ceps \partial_j \ceps,
  \eas
  valid in $\Omega\times (0,\tme)$ for all $i,j\in \{1,2,3\}$,
  we see that
  according to Lemma \ref{lem77}, Lemma \ref{lem80} and Lemma \ref{lem81}, an application of the Aubin-Lions lemma
  (\cite[Ch.~III.2.2]{temam}) provides a sequence $(\eps_j)_{j\in\N} \subset (0,1)$ and limit functions
  $n,c$ and $u$ such that $\eps_j\searrow 0$ as $j\to\infty$ and such that (\ref{82.2})-(\ref{82.7}) hold as well as
  \begin{eqnarray}
	\neps \wto n
	& & 
	\mbox{in } L^\frac{5}{3}_{loc}(\bar\Omega\times [0,\infty)), \label{82.100} \\
	\neps\to n
	& &
	\mbox{in } L^\frac{5}{4}_{loc}(\bar\Omega\times [0,\infty))
	\ \mbox{and a.e.~in } \Omega\times (0,\infty), \label{82.10} \\
	\neps(\cdot,t) \to n(\cdot,t)
	& & \mbox{in } L^\frac{5}{4}(\Omega) \mbox{ for all } t\in (0,\infty)\setminus N, \label{82.103} \\
	\sqrt{\neps} \wto \sqrt{n}
	& & 
	\mbox{in } L^2_{loc}([0,\infty);W^{1,2}(\Omega)), \label{82.106} \\
	\sqrt{\ceps} \to \sqrt{c}
	& &
	\mbox{in } L^2_{loc}([0,\infty);W^{1,2}(\Omega)), \label{82.102} \\
	\sqrt{\ceps(\cdot,t)} \to \sqrt{c(\cdot,t)}
	& & 
	\mbox{in $W^{1,2}(\Omega)$ for all } t\in (0,\infty)\setminus N, \label{82.104} \qquad \mbox{and} \\
	\ueps(\cdot,t) \to u(\cdot,t)
	& & \mbox{in } L^2(\Omega) \mbox{ for all } t\in (0,\infty)\setminus N, \label{82.105}
  \end{eqnarray}
  as $\eps=\eps_j\searrow 0$ with some null set $N\subset (0,\infty)$. 
  Since (\ref{82.100}) entails that for each $T>0$ we have
  \bas
	\int_0^T \io \neps^\frac{5}{3} = \int_0^\infty \io \one_{\Omega\times (0,T)} \neps^\frac{5}{3}
	\to \int_0^\infty \io \one_{\Omega\times (0,T)} n^\frac{5}{3} = \int_0^T \io n^\frac{5}{3}
  \eas
  as $\eps=\eps_j\searrow 0$, it follows that also (\ref{82.1}) holds.\\
  Now in order to verify (\ref{energy}), given a nonnegative $\phi\in C_0^\infty([0,\infty))$ we recall (\ref{776.1}) 
  to obtain
  \bea{e_eps}
	-\int_0^\infty \F[\neps,\ceps,\ueps](t) \cdot \phi'(t)dt
	&+& \frac{1}{K} \int_0^\infty \io \bigg\{\frac{|\nabla\neps|^2}{\neps} + \frac{|\nabla \ceps|^4}{\ceps^3}
	+ |\nabla \ueps|^2 \bigg\}(x,t) \cdot \phi(t) dx dt \nn\\[2mm]
	&\le& \F[n_{0\eps}, c_{0\eps}, u_{0\eps}] \cdot \phi(0)
	+ K\int_0^\infty \phi(t)dt
  \eea
  for all $\eps\in (0,1)$. 
  Here combining (\ref{82.103}), (\ref{82.104}) 
  and (\ref{82.105}) with Lemma \ref{lem70} shows that as $\eps=\eps_j\searrow 0$ we have
  \bea{82.77}
	& & \hspace*{-20mm}
	\io \neps(\cdot,t) \ln \neps(\cdot,) \to \io n(\cdot,t)\ln n(\cdot,t),
	\quad
	\io |\nabla\Psi(\ceps(\cdot,t))|^2 \to \io |\nabla\Psi(c(\cdot,t))|^2
	\quad \mbox{and} \quad \nn\\[2mm]
	& & \hspace*{20mm}
	\ueps(\cdot,t) \to u(\cdot,t)
	\quad \mbox{in } L^2(\Omega)
	\qquad \mbox{for all } t\in (0,\infty)\setminus N
  \eea
  and hence
  \bas
	\F[\neps,\ceps,\ueps](t) \to \F[n,c,u](t)
	\qquad \mbox{for all } t\in (0,\infty)\setminus N,
  \eas
  so that since clearly $\F[\neps,\ceps,\ueps](t) \ge - \frac{|\Omega|}{e}$ for all $\eps\in (0,1)$ and $t>0$ and 
  \bas
	\sup_{\eps\in (0,1)} \sup_{t>0} \, \F[\neps,\ceps,\ueps](t) <\infty
  \eas
  according to Lemma \ref{lem77}, we may invoke the dominated convergence theorem to infer that
  \bas
	\int_0^\infty \F[\neps,\ceps,\ueps](t) \cdot \phi'(t)dt
	\to \int_0^\infty \F[n,c,u](t) \cdot \phi'(t)dt
	\qquad \mbox{as } \eps=\eps_j\searrow 0.
  \eas
  Since (\ref{I1}), (\ref{I2}), Lemma \ref{lem70} and (\ref{I3}) warrant that similarly
  \bas
	\F[n_{0\eps},c_{0\eps},u_{0\eps}]
	\to \F[n_0,c_0,u_0]
	\qquad \mbox{as } \eps\searrow 0,
  \eas
  by nonnegativity of $\phi$ we conclude from (\ref{e_eps}) and the weak convergence
  properties in (\ref{82.106}), (\ref{82.107}) and (\ref{82.7}) that
  \bas
	-\int_0^\infty \F[n,c,u](t) \cdot \phi'(s)dt
	&+& \frac{1}{K} \int_0^\infty \io \bigg\{\frac{|\nabla n|^2}{n} + \frac{|\nabla c|^4}{c^3}
	+ |\nabla u|^2 \bigg\}(x,t) \cdot \phi(t) dx dt \\
	&\le& \F[n_0,c_0,u_0] \cdot \phi(0)
	+ K\int_0^\infty \phi(t)dt
  \eas
  for any such $\phi$, and that hence (\ref{energy}) holds.\abs
  Next,
  to deduce (\ref{energy1}) we integrate (\ref{75.1}) in time to see that	
  \be{82.1000}
	\frac{1}{2} \io |\ueps(\cdot,t)|^2 + \int_{t_0}^t \io |\nabla\ueps|^2
	= \frac{1}{2} \io |\ueps(\cdot,t_0)|^2
	+ \int_{t_0}^t \io \neps \ueps\cdot \nabla \Phi
	\qquad \mbox{for all $t_0\ge 0$ and $t>t_0$,}
  \ee
  where for any such $t_0$ and $t$ we have
  \bas
	\int_{t_0}^t \io \neps \ueps\cdot \nabla \Phi
	\to \int_{t_0}^t \io n u\cdot\nabla \Phi
	\qquad \mbox{as } \eps=\eps_j\searrow 0,
  \eas
  because
  \be{82.188}
	\neps \ueps \to nu
	\qquad \mbox{in } L^1_{loc}(\bar\Omega\times [0,\infty))
	\qquad \mbox{as } \eps=\eps_j\searrow 0
  \ee
  due to (\ref{82.1}), (\ref{82.66}) and the fact that $\frac{3}{5}+\frac{3}{10}=\frac{9}{10}<1$.
  Therefore, (\ref{82.7}) together with the last property in (\ref{82.77}) ensures that if 
  $t_0\in [0,\infty)\setminus N$ and $t>t_0$ then indeed
  \bas
	\frac{1}{2} \io |u(\cdot,t)|^2 + \int_{t_0}^t \io |\nabla u|^2
	&\le& \lim_{\eps=\eps_j\searrow 0} \bigg\{ \frac{1}{2} \io |\ueps(\cdot,t_0)|^2
	+ \int_{t_0}^t \io \neps \ueps\cdot \nabla \Phi \bigg\} \\
	&=& \frac{1}{2} \io |u(\cdot,t_0)|^2
	+ \int_{t_0}^t \io n u\cdot\nabla \Phi,
  \eas
  as desired.\\
  Now in order to verify that $(n,c,u)$ is a global weak solution of 
  (\ref{0}), (\ref{0i}), (\ref{0b}) in the sense of Definition \ref{defi_weak},
  we first note that the regularity properties (\ref{reg_weak}) therein are obvious from (\ref{82.1}), (\ref{82.2}),
  (\ref{82.4}), (\ref{82.107}), (\ref{82.6}) and (\ref{82.7}), that clearly $n$ and $c$ inherit nonnegativity
  from $\neps$ and $\ceps$, and that $\nabla\cdot u=0$ a.e.~in $\Omega\times (0,\infty)$ according to (\ref{0eps}) and
  (\ref{82.7}).
  To prepare a derivation of (\ref{reg_weak2}) and (\ref{w1})-(\ref{w3}), we observe that
  in view of the dominated convergence theorem, 
  \bas
	F_\eps'(\neps) \chi(\ceps) \ceps^\frac{3}{4} \to \chi(c) c^\frac{3}{4}
	\qquad \mbox{in } L^\frac{20}{3}_{loc}(\bar\Omega\times  [0,\infty))
	\qquad \mbox{as } \eps=\eps_j\searrow 0
  \eas
  as a consequence of (\ref{82.1}) and (\ref{82.3}), the boundedness of $(\ceps)_{\eps\in (0,1)}$ in 
  $L^\infty(\Omega\times (0,\infty))$ and the fact that $F'_\eps\nearrow 1$ on $(0,\infty)$ as $\eps\searrow 0$ 
  by (\ref{Feps2}).
  Combining this with (\ref{82.1}) and (\ref{82.107}) shows that
  \be{82.15}
	\neps F_\eps'(\neps) \chi(\ceps) \nabla\ceps
	= 4\neps \cdot F_\eps'(\neps) \chi(\ceps) \ceps^\frac{3}{4} \cdot \nabla\ceps^\frac{1}{4}
	\wto 4n \cdot \chi(c) c^\frac{3}{4} \cdot \nabla c^\frac{1}{4}
	= n\chi(c)\nabla c
	\qquad \mbox{in } L^1_{loc}(\bar\Omega\times [0,\infty))
  \ee
  as $\eps=\eps_j\searrow 0$.
  We proceed to make sure that
  \be{82.16}
	F_\eps(\neps)\to n
	\qquad \mbox{in } L^\frac{5}{3}_{loc}(\bar\Omega\times  [0,\infty))
	\qquad \mbox{as } \eps=\eps_j\searrow 0.
  \ee
  Indeed, for each fixed $T>0$ we have
  \bea{82.166}
	\|F_\eps(\neps)-n\|_{L^\frac{5}{3}(\Omega\times (0,T))}
	&\le& 
	\|F_\eps(\neps)-F_\eps(n)\|_{L^\frac{5}{3}(\Omega\times (0,T))}
	+ \|F_\eps(n)-n\|_{L^\frac{5}{3}(\Omega\times (0,T))} \nn\\
	&\le& \|F_\eps'\|_{L^\infty((0,\infty))} \|\neps-n\|_{L^\frac{5}{3}(\Omega\times (0,T))}
	+ \|F_\eps(n)-n\|_{L^\frac{5}{3}(\Omega\times (0,T))},
  \eea
  where $\|\neps-n\|_{L^\frac{5}{3}(\Omega\times (0,T))} \to 0$ as $\eps=\eps_j\searrow 0$ by (\ref{82.1}),
  and where since
  \bas
	\Big\|F_\eps(n(\cdot,t))-n(\cdot,t)\Big\|_{L^\frac{5}{3}(\Omega)}^\frac{5}{3} 
  	\le 2^\frac{5}{3}\|n(\cdot,t)\|_{L^\frac{5}{3}(\Omega)}^\frac{5}{3}
	\qquad \mbox{ for a.e.~$t>0$}
  \eas
  due to (\ref{Feps1}), (\ref{82.1})
  guarantees that also
  \bas
  	\int_0^T \Big\|F_\eps(n(\cdot,t))-n(\cdot,t)\Big\|_{L^\frac{5}{3}(\Omega)}^\frac{5}{3} dt \to 0
	\qquad \mbox{as } \eps=\eps_j\searrow 0
  \eas
  by the dominated convergence theorem.
  As $0\le F_\eps' \le 1$ by (\ref{Feps1}), (\ref{82.166}) and (\ref{82.1}) prove (\ref{82.16}), 
  from which we particularly obtain that
  \be{82.17}
	F_\eps(\neps) f(\ceps) \to nf(c)
	\qquad \mbox{in } L^1_{loc}(\bar\Omega\times [0,\infty))
	\qquad \mbox{as } \eps=\eps_j\searrow 0,
  \ee
  because again Lebesgue's theorem along with (\ref{82.3}) and (\ref{cinfty_eps}) ensures that
  $f(\ceps) \to f(c)$ in $L^\frac{5}{2}_{loc}(\bar\Omega\times [0,\infty))$ as $\eps=\eps_j\searrow 0$.\\
  Next, since 		
  (\ref{82.3}) and (\ref{cinfty_eps}) furthermore imply that
  $\ceps\to c$ in $L^\frac{10}{7}_{loc}(\bar\Omega\times [0,\infty))$
  as $\eps=\eps_j\searrow 0$, from (\ref{82.66}) we infer that
  \be{82.18}
	\ceps \ueps \to cu
	\qquad \mbox{in } L^1_{loc}(\bar\Omega\times [0,\infty))
	\qquad \mbox{as } \eps=\eps_j\searrow 0.
  \ee
  Finally, following an argument from \cite[Theorem V.3.1.1]{sohr} we use that for each $\varphi\in L^2_\sigma(\Omega)$
  we have $\|\yeps\varphi\|_{L^2(\Omega)} \le \|\varphi\|_{L^2(\Omega)}$ and $\yeps \varphi \to \varphi$ in $L^2(\Omega)$
  as $\eps\searrow 0$ to infer from (\ref{82.77}) that for each $t\in (0,\infty)\setminus N$ we have
  \bas
	\Big\|\yeps \ueps(\cdot,t)-u(\cdot,t) \Big\|_{L^2(\Omega)}
	&\le& \Big\| \yeps\Big(\ueps(\cdot,t)-u(\cdot,t)\Big)\Big\|_{L^2(\Omega)}
	+ \Big\|\yeps u(\cdot,t)-u(\cdot,t) \Big\|_{L^2(\Omega)} \\
	&\le& \Big\|\ueps(\cdot,t)-u(\cdot,t)\Big\|_{L^2(\Omega)}
	+ \Big\|\yeps u(\cdot,t)-u(\cdot,t) \Big\|_{L^2(\Omega)} \\[2mm]
	&\to& 0
	\qquad \mbox{as } \eps=\eps_j\searrow 0,
  \eas
  and that moreover
  \bas
	\Big\|\yeps \ueps(\cdot,t)-u(\cdot,t) \Big\|_{L^2(\Omega)}^2
	&\le& \Big( \|\yeps\ueps(\cdot,t)\|_{L^2(\Omega)} + \|u(\cdot,t)\|_{L^2(\Omega)} \Big)^2 \\
	&\le& \Big(\|\ueps(\cdot,t)\|_{L^2(\Omega)} + \|u(\cdot,t)\|_{L^2(\Omega)} \Big)^2 \\[1mm]
	&\le& 4\sup_{\eps'\in (0,1)} \|u_{\eps'}\|_{L^2(\Omega\times (0,\infty))}^2
	\qquad \mbox{for all $t\in (0,\infty)\setminus N$ and } \eps \in (0,1).
  \eas
  In view of (\ref{77.01}), once more thanks to the dominated convergence theorem this entails that for all $T>0$ we obtain
  \bas
	\int_0^T \Big\|\yeps\ueps(\cdot,t)-u(\cdot,t)\Big\|_{L^2(\Omega)}^2 dt \to 0
	\qquad \mbox{as } \eps=\eps_j\searrow 0.
  \eas
  Thus, 
  \bas
	\yeps\ueps \to u
	\quad \mbox{in } L^2_{loc}(\bar\Omega\times [0,\infty))
	\qquad \mbox{as } \eps=\eps_j\searrow 0,
  \eas
  which in conjunction with (\ref{82.6}) entails that
  \be{82.1888}
	\yeps \ueps \mult \ueps \to u\mult u
	\quad \mbox{in } L^1_{loc}(\bar\Omega\times [0,\infty))
	\qquad \mbox{as } \eps=\eps_j\searrow 0.
  \ee
  Now (\ref{82.188}), (\ref{82.15}), (\ref{82.17}), (\ref{82.18}) and (\ref{82.1888})
  firstly warrant that the integrability requirements in 
  (\ref{reg_weak2}) are satisfied, and secondly, together with (\ref{82.1})-(\ref{82.7}) and (\ref{I1})-(\ref{I3}), 
  allow for passing to the
  limit in the respective weak formulations associated with the equations in (\ref{0eps}).
  In fact, if for $\eps\in (0,1)$
  we multiply the first equation in (\ref{0eps}) by an arbitrary $\phi\in C_0^\infty(\bar\Omega\times [0,\infty))$
  and integrate by parts, then in the resulting identity
  \bas
	-\int_0^\infty \io \neps\phi_t - \io n_{0\eps}\phi(\cdot,0)
	&=& - \int_0^\infty \io \nabla \neps\cdot \nabla \phi
	+ \int_0^\infty \io \neps F_\eps'(\neps)\chi(\ceps) \nabla \ceps \cdot\nabla\phi \\
	& & + \int_0^\infty \io \neps \ueps \cdot \nabla\phi
  \eas
  we may apply (\ref{82.1}), (\ref{I1}), (\ref{82.2}), (\ref{82.15}) and (\ref{82.188}) to take
  $\eps=\eps_j\searrow 0$ in the first, second, third, fourth and fifth integral, respectively, to conclude that
  (\ref{w1}) holds. 
  Likewise, since for all $\phi\in C_0^\infty(\bar\Omega\times [0,\infty))$ and $\eps\in (0,1)$ we have
  \bas
	-\int_0^\infty \io \ceps\phi_t - \io c_{0\eps}\phi(\cdot,0)
	= - \int_0^\infty \io \nabla \ceps\cdot \nabla \phi
	- \int_0^\infty \io F_\eps(\neps)f(\ceps) \phi
	+ \int_0^\infty \io \ceps\ueps\cdot\nabla \phi,
  \eas
  invoking (\ref{82.3}), (\ref{I2}), (\ref{82.107}), (\ref{82.17}) and (\ref{82.18}) and again applying the dominated
  convergence theorem along with (\ref{cinfty_eps}) establishes (\ref{w2}). 
  Finally, 
  given any $\phi\in C_0^\infty(\Omega\times [0,\infty);\R^3)$ satisfying $\nabla\cdot \phi\equiv 0$, from (\ref{0eps})
  we obtain
  \bas
	-\int_0^\infty \io \ueps\cdot\phi_t - \io u_{0\eps}\cdot \phi(\cdot,0)
	= - \int_0^\infty \io \nabla \ueps \cdot\nabla \phi
	+ \int_0^\infty \io \yeps \ueps \mult \ueps \cdot \nabla \phi
	+ \int_0^\infty \io \neps\nabla \Phi\cdot \phi
  \eas
  for all $\eps\in (0,1)$, so that taking $\eps=\eps_j\searrow 0$ and using (\ref{82.6}), (\ref{I3}),
  (\ref{82.7}), (\ref{82.1888}) and (\ref{82.1}) yields (\ref{w3}) and thereby completes the proof.
\qed
\proofc of Theorem \ref{theo_global}. \quad
  The statement is evidently implied by Lemma \ref{lem82}.
\qed

\end{document}